\documentclass[12pt]{article}

\usepackage{amssymb, amsmath}
\usepackage[normalem]{ulem}
\usepackage{amsfonts}
\usepackage{graphicx}
\usepackage{setspace} 
\usepackage{rotating}
\usepackage{enumitem}
\allowdisplaybreaks
\usepackage{hyperref}
\renewcommand{\theequation}{\thesection .\arabic{equation}}
\usepackage{tikz}
\usepackage{mathrsfs}
\usepackage{url}
\usepackage[authoryear]{natbib}
\newcommand{\be}{\begin{equation}}
	\newcommand{\ee}{\end{equation}}
\newcommand{\bea}{\begin{eqnarray}}
	\newcommand{\eea}{\end{eqnarray}}
\newcommand{\bean}{\begin{eqnarray*}}
	\newcommand{\eean}{\end{eqnarray*}}
\newcommand{\brray}{\begin{array}}
	\newcommand{\erray}{\end{array}}
\newcommand{\ben}{\begin{equation}{nonumber}}
	\newcommand{\een}{\end{equation}{nonumber}}

\newtheorem{dfn}{Definition}[section]
\newtheorem{thm}[dfn]{Theorem}
\newtheorem{lmma}[dfn]{Lemma}
\newtheorem{ppsn}[dfn]{Proposition}
\newtheorem{crlre}[dfn]{Corollary}
\newtheorem{xmpl}[dfn]{Example}
\newtheorem{rmrk}[dfn]{Remark}

\newcommand{\bdfn}{\begin{dfn}}
	\newcommand{\bthm}{\begin{thm}}
		\newcommand{\blmma}{\begin{lmma}}
			\newcommand{\bppsn}{\begin{ppsn}}
				\newcommand{\bcrlre}{\begin{crlre}}
					\newcommand{\bxmpl}{\begin{xmpl}}
						\newcommand{\brmrk}{\begin{rmrk}}
							\newcommand{\edfn}{\end{dfn}}
						\newcommand{\ethm}{\end{thm}}
					\newcommand{\elmma}{\end{lmma}}
				\newcommand{\eppsn}{\end{ppsn}}
			\newcommand{\ecrlre}{\end{crlre}}
		\newcommand{\exmpl}{\end{xmpl}}
	\newcommand{\ermrk}{\end{rmrk}}
\newcommand{\p}{\mathbb{P}}

\newcommand{\IR}{\mathbb{R}}

\newcommand{\inpr}[2]{\left\langle#1 \,,\, #2\right\rangle_{\cal{H}}}
\newcommand{\R}{\mathbb{R}}
\newtheorem{theorem}{Theorem}[section]

\newtheorem{lemma}{Lemma}[section]

\newtheorem{proposition}{Proposition}[section]

\def\a*{{\cal A}_{h,*}}
\def\B{{\cal B}(h)}
\def\B1{{\cal B}_1(h)}
\def\b{{\cal B}^{\rm s.a.}(h)}
\def\b1{{\cal B}^{\rm s.a.}_1(h)}

\def \qed {$\Box$}

\usepackage{geometry}
\usepackage{xr}
\makeatletter
\newcommand*{\addFileDependency}[1]{
  \typeout{(#1)}
  \@addtofilelist{#1}
  \IfFileExists{#1}{}{\typeout{No file #1.}}
}
\makeatother

\newcommand*{\myexternaldocument}[1]{
    \externaldocument{#1}
    \addFileDependency{#1.tex}
    \addFileDependency{#1.aux}
}
\myexternaldocument{Main_Bernoulli}
\begin{document}




  \title{Testing Independence of Infinite Dimensional Random Elements: A Sup-norm Approach}
\author{\small 
	Suprio Bhar \\
	\small IIT Kanpur\\
	\small Department of Mathematics and Statistics \\
	\small  Kanpur 208016, India\\
	{\small email: suprio@iitk.ac.in}\\
	\and
	\small Subhra Sankar Dhar \\
	\small  IIT Kanpur\\
	\small   Department of Mathematics and Statistics \\
	\small Kanpur 208106, India\\
	{\small email: subhra@iitk.ac.in}\\
}
\date{}
\maketitle

\setcounter{equation}{0}

\setcounter{section}{0}\ \\


\begin{abstract}
In this article, we study the test for independence of two random elements $X$ and $Y$ lying in an infinite dimensional space ${\cal{H}}$ (specifically, a real separable Hilbert space equipped with the inner product $\langle ., .\rangle_{\cal{H}}$). In the course of this study, a measure of association is proposed based on the sup-norm difference between the joint probability density function of the bivariate random vector $(\langle l_{1}, X \rangle_{\cal{H}}, \langle l_{2}, Y \rangle_{\cal{H}})$  and the product of marginal probability density functions of the random variables $\langle l_{1}, X \rangle_{\cal{H}}$ and $\langle l_{2}, Y \rangle_{\cal{H}}$, where $l_{1}\in{\cal{H}}$ and $l_{2}\in{\cal{H}}$ are two arbitrary elements. It is established that the proposed measure of association equals zero if and only if the random elements are independent. In order to carry out the test whether $X$ and $Y$ are independent or not, the sample version of the proposed measure of association is considered as the test statistic after appropriate normalization, and the asymptotic distributions of the test statistic under the null and the local alternatives are derived. The performance of the new test is investigated for simulated data sets and the practicability of the test is shown for three real data sets related to climatology, biological science and chemical science. 
\end{abstract}


\section{Introduction}\label{Introduction}
\subsection{Key Ideas and Literature Review}
For univariate and multivariate data, there have been several attempts to test whether two or more random variables or vectors are independent or not in various situations (see, e.g., \cite{Blum1961}, \cite{Szekely2007}, \cite{Genest2007}, \cite{Einmahl2008}, \cite{Dette2013},  \cite{Bergsma2014}, \cite{DharEJS2016}, \cite{Han2017}, \cite{Dhar2018}, \cite{Drton2020}, \cite{Chatterjee2021}, \cite{Berrett2021}, \cite{Shi2022}, \cite{Shi2022JASA} and a few references therein). To the best of our knowledge, all these tests are based on a certain measure of association, which equals zero if and only if the finite dimensional random variables or vectors are independent, and consequently, the tests based on those measures of associations lead to consistent tests. Therefore, one first needs to develop a measure of association having such {\it necessary and sufficient} relation with independence in the infinite dimensional setting. Once such a relation is found, in principle, it allows us to formulate a reasonable test statistic for checking whether two random elements are independent or not in infinite dimensional space. This article addresses this issue in the subsequent sections.

Let us first discuss the utility of the framework in terms of infinite dimensional spaces. In many interdisciplinary subjects like Biology, Economics, Climatology or Finance, we recently come across many data sets, where the dimension of the data is much larger than the sample size of the data set, and in majority of the cases, the standard multivariate techniques cannot be implemented because of their high dimensionalities relative to the sample sizes. In order to overcome this problem, one can embed such data into a suitable infinite dimensional space. For instance, functional data (see \cite{Ramsey2002}, \cite{Fer2006}) is an infinite dimensional data, and one can analyse the functional data using the techniques adopted for infinite dimensional data. In this work, we consider random elements lying in a real separable Hilbert space. The separability of the space allows us to write the random elements as linear combinations of a suitable orthonormal basis, with random coefficients, making it easier to handle the theoretical issues. In the context of measure of association or the test for independence of infinite dimensional random elements, we would like to mention that only a limited number of articles are available in the literature on this topic, and the contributions of the major ones are described here.

Almost a decade ago, \cite{Russell2013} (see also \cite{Russell2018} and \cite{Russell2021}) studied the distance covariance, which is proposed by \cite{Szekely2007} for finite dimensional random vectors, in a certain metric space. However, \cite{Russell2013} did not study the corresponding testing of hypothesis problems based on the proposed measure. Before the aforesaid work, \cite{Juan2006} studied the random projection based goodness of fit test and two sample tests for the infinite dimensional element, which can be applied on the test for  independence of the infinite dimensional random elements as well. During the almost same period, \cite{Gretton2007} proposed a new methodology based on the Hilbert-Schmidt independence criterion (HSIC) for testing independence of two finite dimensional random vectors, and it can too be applied for the Hilbert space valued random elements. However, none of these articles used the density functions of the random variables obtained by the projections of the random elements lying in an infinite dimensional space. This article proposes a methodology for testing the independence of random elements defined in an infinite dimensional space, specifically a real separable Hilbert space, based on the probability density functions of projections of said random elements. 

\subsection{Contributions}
We consider two random elements $X$ and $Y$ lying in a real separable Hilbert space ${\cal{H}}$ (equipped with the inner product $\langle ., . \rangle_{\cal{H}}$) and propose a methodology of testing independence between $X$ and $Y$. This is the first major contribution of this article. The methodology is developed based on sup-norm distance between the joint probability density of the bivariate random vector $(\langle l_{1}, X \rangle_{\cal{H}}, \langle l_{2}, Y \rangle_{\cal{H}})$  and the product of marginal probability density functions of the random variables $\langle l_{1}, X \rangle_{\cal{H}}$ and $\langle l_{2}, Y \rangle_{\cal{H}}$, where $l_{1}\in{\cal{H}}$ and $l_{2}\in{\cal{H}}$ are two arbitrary non-zero elements. This measure of association equals with zero for all $l_{1}\in{\cal{H}}$ and $l_{2}\in{\cal{H}}$ {\it if and only if} $X$ and $Y$ are independent random elements, and it is non-negative as well. 

The second major contribution is the following. For a given data, a sample version of the measure with appropriate normalization is proposed, which is considered as the test statistic for testing independence of infinite dimensional random elements,  and the asymptotic distribution of the sample version under null and local/contiguous alternatives is derived after appropriate normalization. These theoretical results enable us to study the consistency and the asymptotic power of the corresponding test. Unlike other tests, since the test statistic used in this methodology depends on the choice of the kernel and the associated bandwidth, the performance of the proposed test can be enhanced by suitable choices of the kernel and the bandwidth. 

The third major contribution of the work is the implementation of the test. As the test statistic is based on the supremum over the infinite choices of infinite dimensional random elements $l_1$ and $l_2$, the exact computation of the test statistic is intractable for a given data. To overcome this problem, $l_1$ and $l_2$ are chosen over certain finite collection of possibilities of the choices, and it is shown that this method approximates the actual statistic when the number of possible choices of $l_1$ and $l_2$ are sufficiently large. Overall, it is established that the approximated test is easily implementable, and it gives satisfactory results in analysing real data as well. 

\subsection{Challenges}
In the course of this work, we overcome a few mathematical challenges. The first challenge involves the interpretation of independence between two infinite dimensional random elements, and it is resolved using the concept of projection towards all possible directions (see Lemma \ref{lemindep}). Using Proposition \ref{equivalence}, we further reduce the problem by relying on vectors of unit length, rather that those with arbitrary lengths. This procedure reduces the complexity of the optimization problem associated with infinite dimensional projection vector to a large extent. The second challenge concerns the optimization problem associated with time parameters involved with the test statistic. To overcome this issue, we approximate the test statistic in two steps: first, we use a finite dimensional approximation for any random vector on an orthonormal basis and then, we compute the test statistic on a sufficiently large number of time points chosen over a sufficiently large interval. Using advanced techniques of analysis, it is shown that the later version of the test statistic approximates arbitrarily well the original test statistic (see Proposition \ref{GL} and Lemma \ref{from-max-on-grid-to-sup}). The third challenge is related to asymptotic distribution of the test statistic. As the key term of the test statistic involves a triangular array and a handful number of algebraically complex terms, the asymptotic distribution of the key term follows after careful use of CLT associated with a triangular array (see the proofs of Lemmas \ref{ASI}, \ref{ASII}, \ref{ASIII}, \ref{pointwiseasymptotic} and \ref{pointwiseasymptotic-multidim}). Finally, the asymptotic distribution of the test statistic follows from Cramér–Wold theorem and continuous mapping theorem (see \cite{van1998}).

\subsection{Applications}
We analyse well known climate data, which consists of daily temperature and precipitation at thirty five locations in Canada averaged over 1960 to 1994 (see \cite{Ramsey2005}). Note that the dimension of the data equals 365, which is much larger than the sample size (i.e., 35) of the data, and embedding such a high dimensional data into a specific Hilbert space is legitimate enough. From the point of view of climate science, it is of interest to know how far the temperature depends on the precipitation at a given time point. The analysis of this real data using the proposed methodology gives some insight about the long standing issue in climatology.

Along with it, two more data sets are also analysed to show the practicability of the proposed methodology. One data set is the well-known Berkeley growth data, which consists the heights of 39 boys and 54 girls measured at 31 time points between the ages 1 and 18 years, and from medical science's point of view, it is of interest to know whether the growth (in terms of height) of the boys and the girls are independent or not. We address it using our proposed test. Another data set is Coffee data, which contains spectroscopy readings taken
at 286 wavelength values for 14 samples of each of the two varieties of coffee, namely,
Arabica and Robusta. As coffee is one of the most popular beverages, one may be interested to know whether the chemical features of the Arabica coffee and the Robusta coffee are independent or not, and this issue is also addressed using our test.

\subsection{Organization of the Article}
The article is organized as follows. Section \ref{Methodology} introduces the basic concepts of separable Hilbert space valued random elements, and Section \ref{PCI} proposes the criterion for checking independence between two such random elements. Afterwards, for a given data, the estimation of the criterion is studied in Section \ref{ET}, and Section \ref{FTS} formulates the test statistic. The asymptotic properties of the test statistic and associated test are investigated in Section \ref{APT}, and Section \ref{APS} studies the performance of the proposed test in terms of the asymptotic power for various examples. The performance of the proposed test for finite sample sizes is studied in Section \ref{FSLPS}, and real data analysis is carried out in Section \ref{RDA}. Section \ref{CR} contains some concluding remarks, and finally, Sections \ref{TD} and \ref{ANS} contain all technical details and additional numerical results, respectively.

\section{Methodology}\label{Methodology}
We first formulate the statement of the problem and subsequently construct an appropriate test statistics. Let ${\cal{H}}$ denote a real separable Hilbert space equipped with the inner product $\langle ., . \rangle_{\cal{H}}$ and the norm $||.||_{\cal{H}} = \sqrt{\langle ., . \rangle_{\cal{H}}}$, and suppose that $(\Omega, {\cal{A}}, P)$ is a probability space. 
Let us now consider an ${\cal{H}}\times{\cal{H}}$ valued bivariate random element $(X, Y)$, which is a measurable mapping from $(\Omega, {\cal{A}}, P)$ into ${\cal{H}}\times{\cal{H}}$ equipped with its Borel $\sigma$-algebra ${\cal{B}}({\cal{H}}\times{\cal{H}})$ generated by the open sets of ${\cal{H}}\times{\cal{H}}$. In other words, for any Borel set $U_{1}\in{\cal{H}}$ and $U_{2}\in{\cal{H}}$, $(X^{-1}(U_{1}), Y^{-1}(U_{2}))\in{\cal{A}}$. We now want to check whether the random elements $X$ and $Y$ are independent or not, and in this regard, a new criterion for checking independence is proposed here. 

\subsection{Proposed Criterion for Independence}\label{PCI}
We here propose a criteria for checking independence based on certain one-dimensional projections of the random elements $X$ and $Y$ lying on a separable Hilbert space ${\cal{H}}$. In order to formulate the criterion, a useful lemma is stated here. 

\begin{lemma}\label{lemindep}
Let $X$ and $Y$ be two ${\cal{H}}$ valued random elements on a probability space $(\Omega, \cal{A}, P)$, where ${\cal{H}}$ is a real separable Hilbert space. Then ${\cal{H}}$-valued $X$ and ${\cal{H}}$-valued $Y$ are two independent random elements if and only if the real valued random variable $\langle l_{1}, X \rangle_{{\cal{H}}}$ and the real valued random variable $\langle l_{2}, Y \rangle_{{\cal{H}}}$ are independent for every $l_{1}\in {\cal{H}}$ and $l_{2}\in{\cal{H}}$. In other words, 
\[P(\inpr{l_1}{X} \in A, \inpr{l_2}{Y} \in B) = P(\inpr{l_1}{X} \in A) P(\inpr{l_2}{Y} \in B)\]
for all Borel subsets $A$ and $B$ of $\R$ if and only if
\[P(X \in C, Y \in D) = P(X \in C) P(Y \in D)\]
for all Borel subsets $C$ and $D$ of $\cal{H}$.
\end{lemma}

The assertion in Lemma \ref{lemindep} indicates that $P_{X, Y}(C \times D) = P_{X} (C) P_{Y} (D)$ for all Borel sets $C\in{\cal{H}}$ and $D\in {\cal{H}}$, where $P_{X}$ and $P_{Y}$ are marginal measures associated with $X$ and $Y$, respectively of the joint law $P_{X, Y}$ if and only if $Q(A\times B) = Q_{X}(A) Q_{Y}(B)$, where $Q$ is the joint law of ($\langle l_{1}, X \rangle_{{\cal{H}}} , \langle l_{2}, Y \rangle_{{\cal{H}}}$), $Q_{X}$ is the measure associated with $\langle l_{1}, X \rangle_{{\cal{H}}}$, and $Q_{Y}$ is the measure associated with $\langle l_{2}, Y \rangle_{{\cal{H}}}$, respectively. Here $A\in {\cal{B}}(\mathbb R)$ and $B\in{\cal{B}}(\mathbb R)$ are two arbitrary Borel sets. This fact indicates that if $\langle l_{1}, X \rangle_{{\cal{H}}}$ and $\langle l_{2}, Y \rangle_{{\cal{H}}}$ are absolutely continuous random variables, then the joint probability density function of $\langle l_{1}, X \rangle_{{\cal{H}}}$ and $\langle l_{2}, Y \rangle_{{\cal{H}}}$ equals with the product of marginal probability density functions of $\langle l_{1}, X \rangle_{{\cal{H}}}$ and $\langle l_{2}, Y \rangle_{{\cal{H}}}$ for all $l_1$ and $l_2$. Writing with notation, suppose that $f_{(\langle l_{1}, X \rangle_{{\cal{H}}}, \langle l_{2}, Y \rangle_{{\cal{H}}})}(., .)$, $f_{\langle l_{1}, X \rangle_{{\cal{H}}}} (.)$ and $f_{\langle l_{2}, Y \rangle_{{\cal{H}}}}(.)$ are the joint probability density, marginal probability density functions of  ($\langle l_{1}, X \rangle_{{\cal{H}}} , \langle l_{2}, Y \rangle_{{\cal{H}}}$), $\langle l_{1}, X \rangle_{{\cal{H}}}$ and $\langle l_{2}, Y \rangle_{{\cal{H}}}$, respectively. Then $X$ and $Y$ are independent random elements if and only if
\begin{equation}\label{identity}
 f_{(\langle l_{1}, X \rangle_{{\cal{H}}}, \langle l_{2}, Y \rangle_{{\cal{H}}})}(s, t) - f_{\langle l_{1}, X \rangle_{{\cal{H}}}} (s) f_{\langle l_{2}, Y \rangle_{{\cal{H}}}} (t) = 0   
\end{equation} for all $s\in\mathbb{R}$, $t\in\mathbb{R}$, $l_{1}\in{\cal{H}}$ and $l_{2}\in {\cal{H}}$. The identity in \eqref{identity} motivates us to formulate  the following criterion for checking independence between $X$ and $Y$: 
\begin{equation}\label{prilimarycriterion}
T(R_1, R_2):=\sup_{||l_1||_{{\cal{H}}} \leq R_{1}, ||l_2||_{{\cal{H}}} \leq R_{2}}\sup_{s\in\mathbb{R}, t\in\mathbb{R}}|f_{(\langle l_{1}, X \rangle_{{\cal{H}}}, \langle l_{2}, Y \rangle_{{\cal{H}}})}(s, t) - f_{\langle l_{1}, X \rangle_{{\cal{H}}}} (s) f_{\langle l_{2}, Y \rangle_{{\cal{H}}}} (t)|,
\end{equation} where $R_1$ and $R_2$ are two arbitrary large positive real numbers. Using \eqref{identity}, we have the following:

\begin{proposition}
\label{prilimarycriterionstatement}
$T(R_1, R_2) = 0$ if and only if $X$ and $Y$ are $\cal{H}$-valued independent random elements, where $T(R_1, R_2)$ is as defined in \eqref{prilimarycriterion}. {\color{black} In other words, $T(R_1, R_2)$ is degenerate at $0$ if and only if $X$ and $Y$ are $\cal{H}$-valued independent random elements.}
\end{proposition}

However, note that finding the supremum over $\{l_1 : ||l_1||_{\cal{H}} \leq R_1\}$ and $\{l_2 : ||l_2||_{\cal{H}} \leq R_2\}$ is not easily tractable in practice, and to overcome it, one can find the supremum over the boundary of the unit sphere, i.e., $\{l_1 : ||l_1||_{\cal{H}} = 1\}$ and $\{l_2 : ||l_2||_{\cal{H}} = 1\}$ defined in ${\cal{H}}$, and this equivalence is asserted in Proposition \ref{equivalence}.

\begin{proposition}\label{equivalence}
Let 
\begin{equation}\label{finalycriterion}
T:=\sup_{||l_1||_{{\cal{H}}} = 1, ||l_2||_{{\cal{H}}} = 1}\sup_{s\in\mathbb{R}, t\in\mathbb{R}}|f_{(\langle l_{1}, X \rangle_{{\cal{H}}}, \langle l_{2}, Y \rangle_{{\cal{H}}})}(s, t) - f_{\langle l_{1}, X \rangle_{{\cal{H}}}} (s) f_{\langle l_{2}, Y \rangle_{{\cal{H}}}} (t)|.
\end{equation} Then $$T(R_1, R_2) = 0\Leftrightarrow T = 0,$$ where $T(R_1, R_2)$ is the same as defined in \eqref{prilimarycriterion}.
\end{proposition}
Finally, the following theorem characterizes the independence of $X$ and $Y$ based on equality of $T$ with zero. 

\begin{theorem}\label{characterization}
Let $(X, Y)$ be a bivariate random element taking values on ${\cal{H}}\times {\cal{H}}$, where ${\cal{H}}$ is a real separable Hilbert space. Then, $X$ and $Y$ are independent if and only if $T = 0$, where $T$ is defined in \eqref{finalycriterion}.
\end{theorem}

\begin{rmrk}
In our framework, we consider a bivariate random element $(X, Y)$ taking values on ${\cal{H}}\times {\cal{H}}$ such that the real valued random variables $\langle l_{1}, X \rangle_{{\cal{H}}}$ and $\langle l_{2}, Y \rangle_{{\cal{H}}}$ are absolutely continuous, for all $l_1$ and $l_2$ in $\cal H$. It may also be viable to work when the said real valued random variables are discrete.
\end{rmrk}

The assertion in Theorem \ref{characterization} implies that the testing of independence between the random elements $X$ and $Y$ defined in ${\cal{H}}$ is equivalent to testing $T = 0$, and hence, in order to test whether $T = 0$ or not (or/equivalently $X$ and $Y$ are independent random elements or not) for a given data, one needs to have an {\it appropriate} estimator of $T$.

\subsection{Estimation of $T$}\label{ET}
Let $(X_1, Y_1), \ldots, (X_n, Y_n)$ be an i.i.d sequence of random elements, and they are identically distributed with $(X, Y)$. In order to estimate $T$, one first needs to find $l_1$ and $l_2$ from the unit sphere in ${\cal{H}}$;  see \eqref{finalycriterion}. However, in \eqref{finalycriterion}, $l_1$ and $l_2$ appear in the joint and marginal probability density functions and as such, we end up with a two fold issue. First one involves an infinite-dimensional optimization in $l_1$ and $l_2$ and the second one is the estimation of the corresponding joint and marginal probability density functions. The following may be one of the procedures.  

\subsubsection{Approximation of $l_1$ and $l_2$}\label{approxl}
In view of the fact that ${\cal{H}}$ is a real separable Hilbert space, and since $l_{1}\in{\cal{H}}$ and $l_{2}\in{\cal{H}}$, we have 
$$l_{1} = \sum\limits_{i = 1}^{\infty}l_{1, i} e_{i} ~\mbox{and}~l_{2} = \sum\limits_{i = 1}^{\infty} l_{2, i} e_{i}, $$ where for a fixed $i$, $e_{i} = (\underbrace {0, \ldots, 0}_{i - 1}, 1, 0, \ldots)$ is an infinite dimensional basis in ${\cal{H}}$, and $l_{1, i}$ and $l_{2, i}$ ($i = 1, 2, \ldots$) are the coefficients of the orthonormal expansion of $l_1$ and $l_2$. Note that
\[\sum\limits_{i = 1}^{\infty} l_{1, i}^{2} = 1, \text{ and } \sum\limits_{i = 1}^{\infty} l_{2, i}^{2} = 1,\]
and
$$\left|\left|l_{1} - \sum\limits_{i = 1}^{M} l_{1, i} e_{i}\right|\right|^{2}_{\cal{H}}\rightarrow 0 ~\mbox{and}~ \left|\left|l_{2} - \sum\limits_{i = 1}^{M} l_{2, i} e_{i}\right|\right|^{2}_{\cal{H}}\rightarrow 0$$ as $M\rightarrow\infty$ as ${\cal{H}}$ is a real separable Hilbert space (see, e.g., \cite{Rudin1991}). Using this fact, we propose an approximate choices of $l_{j}$ ($j = 1$ and $2$) as 
\begin{equation}\label{approximatel}
\hat{l}_{j}^{M} = \sum\limits_{i = 1}^{M} l_{j, i} e_{i},
\end{equation} where $M$ is a sufficiently large positive integer. Afterwards, to compute the criterion $T$ or estimate $T$, we compute the supremum over $(l_{j, 1}, \ldots, l_{j, M})$ for $j = 1$ and 2, where $\sum\limits_{i = 1}^{M} l_{j, i}^{2} = 1$, and this optimization problem can approximately be solved using polar transformation. 

For $j = 1$ and $2$, let us consider the following transformation into the spherical coordinate system.
\begin{equation}\label{lM}
\begin{split}
l_{j, 1} =& \cos\theta_{j, 1},\\
l_{j, 2} =& \sin\theta_{j, 1}\cos\theta_{j, 2},\\
l_{j, 3} =& \sin\theta_{j, 1}\sin\theta_{j, 2}\cos\theta_{j, 3},\\
&\vdots\\
l_{j, M - 1} =& \sin\theta_{j, 1}\sin\theta_{j, 2}\ldots\sin\theta_{j, M - 2}\cos\theta_{j, M - 1},\\
l_{j, M} =& \sin\theta_{j, 1}\sin\theta_{j, 2}\ldots\sin\theta_{j, M - 2}\sin\theta_{j, M - 1},
\end{split}
\end{equation}
where for $i = 1,\ldots, M - 2$, $\theta_{j, i}\in (-\frac{\pi}{2}, \frac{\pi}{2})$  and $\theta_{j, M - 1}\in (- \pi, \pi)$. Observe that for $j = 1$ and 2, $\sum\limits_{i = 1}^{M} l_{j, i}^{2} = 1$. In practice, in the course of implementing the test, we consider $K$ equally spaced choices of $\theta_{j, i}$ from $(-\frac{\pi}{2}, \frac{\pi}{2})$ for $i = 1, \ldots, M-2$, and similarly, $K$ equally spaced $\theta_{j, M - 1}$ are chosen from $(-\pi, \pi)$. Using this grid searching of $\theta_{j, i}$ ($j = 1$ and 2, and $i = 1, \ldots, M - 1$), suppose that the supremum associated with $l_1$ and $l_2$ involved in $T$ (see \eqref{finalycriterion}) attains at some $\theta_{j. i}^{(K)}$ for $j = 1$ and 2, and $i = 1, \ldots, M$, and corresponding associated $l_{j, i}$ are denoted by $l_{j, i}^{(K)}$. Finally, our approximation of $l_{j}$,  depending on $K$, is considered as
\begin{equation}\label{finall}
\hat{l}_{j}^{K} = \sum\limits_{i = 1}^{M} l_{j, i}^{(K)} e_{i}.
\end{equation} Observe that $\sum\limits_{i = 1}^{M} (l_{j, i}^{(K)})^{2} = 1$ for $j = 1$ and 2. The following proposition guarantees that $\hat{l}_{j}^{K}$ is approximate $l_{j}$ ($j = 1$ and 2) arbitrarily well. {\it Note that $\hat{l}_{j}^{K}$ ($j = 1$ and 2) depend on $M$ as well. However, for notational convenience, we do not explicitly write $M$ in the expression of $\hat{l}_{j}^{K}$.}

\begin{proposition}\label{lk}
For $j = 1$ and 2, $\left|\left|\hat{l}_{j}^{K} - l_{j}\right|\right|_{\cal{H}}\rightarrow 0$ as $M\rightarrow\infty$ and $K\rightarrow\infty$, where $\hat{l}_{j}^{K}$ is defined in \eqref{finall}.
\end{proposition}

Observe that $\hat{l}_{j}^{K}$ is also unknown in practice as the joint and marginal probability density functions involved in $T$ (see \eqref{finalycriterion}) are unknown. Therefore, as indicated from Proposition \ref{lk} that $\hat{l}_{j}^{K}$ is a {\it good} approximation of $l_{j}$ ($j = 1$ and 2), and since computing $\hat{l}_{j}^{K}$ is much less complex relative to computing $l_j$ ($j = 1$ and 2), we will work on $\hat{l}_{j}^{K}$ in formulating the test statistic and study its properties.

\subsubsection{Formulation of Test Statistic}\label{FTS}
As we said earlier, formally speaking we want to test $H_{0} : X\perp Y$ against $H_{1} : X\not\perp Y$, where $X$ and $Y$ are two random elements lying in a real separable Hilbert space ${\cal{H}}$, and it follows from Theorem \ref{characterization} that it is equivalent to test $H_{0}^{*} : T = 0$ against $H_{1}^{*}: T\neq 0$, where $T$ is defined in \eqref{finalycriterion}. In order to test $H_{0}^{*}$ against $H_{1}^{*}$, one needs to have an {\it appropriate} estimator of $T$, and as an initial step in the course of this work, Section \ref{approxl} studied how to approximate $l_1$ and $l_2$ involved in $T$ (see \eqref{finalycriterion}), and the approximated $l_j$ ($j = 1$ and 2) are denoted by $\hat{l}_{j}^{K}$ defined in \eqref{finall}. 

Here one first needs to project the infinite dimensional data into one dimensional space through $\hat{l}_1^{K}$ and $\hat{l}_2^{K}$. Let us first consider the data $(X_1, Y_1), \ldots, (X_n, Y_n)$, which belong to the same probability space of $(X, Y)$, and suppose that the transformed data are 
$(\langle \hat{l}_{1}^{K}, X_{1}\rangle_{\cal{H}}, \langle \hat{l}_{2}^{K}, Y_{1}\rangle_{\cal{H}}), \ldots, (\langle \hat{l}_{1}^{K}, X_{n}\rangle_{\cal{H}}, \langle \hat{l}_{2}^{K}, Y_{n}\rangle_{\cal{H}})$, where the projection vectors $\hat{l}_{1}^{K}$ and $\hat{l}_{2}^{K}$ are unknown in practice as discussed at the end of Section \ref{approxl}. In order to identify the optimum $\hat{l}_{1}^{K}$ and $\hat{l}_{2}^{K}$ for a given data, we need to estimate the joint probability density function $f_{(\langle l_{1}, X \rangle_{{\cal{H}}}, \langle l_{2}, Y \rangle_{{\cal{H}}})} (., .)$ and the marginal probability density functions  $f_{\langle l_{1}, X \rangle_{{\cal{H}}}}(.)$ and $f_{\langle l_{2}, Y \rangle_{{\cal{H}}}}(.)$. We describe the estimation of the probability density functions in the following.  

Let $k : \mathbb{R}^{2}\rightarrow (0, \infty)$ be a bivariate function such that $\int k(x, y) dx dy = 1$, $k_1 : \mathbb{R}\rightarrow (0, \infty)$ and $k_2 : \mathbb{R}\rightarrow (0, \infty)$ are such that $\int k_{1} (x) dx = 1$ and $\int k_{2} (x) dx = 1$. Note that $k$, $k_1$ (and $k_2$) can be considered as the kernels involved in the estimator of the bivariate and univariate probability density functions, respectively. For details on kernel density estimation, the readers may refer to \cite{Silverman1998}, and we shall impose a few more technical conditions on the kernels, to be described at the appropriate places. Now, suppose that $T_{n}$ is an appropriate estimator of $T$, which can be formulated as follows. 

\begin{equation}\label{statistic}
\begin{aligned}
T_{n} &= 
\sup_{\substack{\theta_{j, i}\in (-\frac{\pi}{2}, \frac{\pi}{2})\\ j = 1, 2; i = 1, \ldots, M - 1\\ \theta_{j, M - 2}\in (-\pi, \pi)\\ j = 1, 2\\ s\in\mathbb{R}, t\in\mathbb{R}}}
\left|\frac{1}{n h_{n}^2}\sum\limits_{i = 1}^{n}k\left(\frac{\langle \hat{l}_{1}^{K}, X_{i}\rangle_{\cal{H}} - s}{h_{n}}, \frac{\langle \hat{l}_{2}^{K}, Y_{i}\rangle_{\cal{H}} - t}{h_{n}}\right)\right.\\ - &\left.\frac{1}{n^{2}h_{n}^{2}}\sum\limits_{i, j = 1}^{n}k_{1}\left(\frac{\langle \hat{l}_{1}^{K}, X_{i}\rangle_{\cal{H}} - s}{h_{n}}\right)k_{2}\left(\frac{\langle \hat{l}_{2}^{K}, Y_{j}\rangle_{\cal{H}} - t}{h_{n}}\right)\right|, 
\end{aligned}
\end{equation} 
where $h_{n}$ is a sequence of bandwidth such that $h_{n}\rightarrow 0$ as $n\rightarrow\infty$.
Note that $\hat{l}_{j}^{K}$ involves $\theta_{j, i}$ (see \eqref{lM} and \eqref{finall}), and for this reason, supremum taken over $\hat{l}_{j}^{K}$ eventually boils down to the supremum taken over $\theta_{j, i}$ ($j = 1$ and 2, $i = 1, \ldots, M$) in the expression of $T_{n}$ (see \eqref{statistic}). The asymptotic properties of $T_{n}$ are studied in the subsequent section. 

\subsection{Asymptotic Properties of $T_{n}$}\label{APT}
Recall that the testing of hypothesis problem $H_{0} : X\perp Y$ against $H_{1}: X\not\perp Y$ is equivalent to testing of hypothesis problem $H_{0}^{*} : T = 0$ against $H_{1}^{*} : T\neq 0$. In order to check whether for a given data $T = 0$ or not, the test statistic $T_{n}$ (see \eqref{statistic}) is formulated in Section \ref{FTS}, and in order to carry out the test based on $T_{n}$, the distributional feature of $T_{n}$ is required. However, due to complex formulation of $T_{n}$, the derivation of exact distribution of $T_{n}$ may not be tractable, and to overcome it, we here derive the asymptotic distribution of $T_{n}$. We first assume the following technical conditions and then state the asymptotic distribution of $T_{n}$ under local alternatives, and the asymptotic distribution under null hypothesis (i.e., $H_{0}$ or $H_{0}^{*}$).

\noindent (A1) The sequence of bandwidth $\{h_{n}\}_{n\geq 1}$ is such that $h_{n}\rightarrow 0$ as $n\rightarrow\infty$, $nh_{n}^{4}\rightarrow c$ (for some constant $c > 0$) as $n\rightarrow\infty$.

\noindent (A2) The partial derivatives of the joint probability density function $f$ of the bivariate random vector $(\langle l_{1}, X \rangle_{\cal{H}}, \langle l_{2}, Y \rangle_{\cal{H}})$ exist, i.e., $\frac{\partial f_{(\langle l_{1}, X \rangle_{\cal{H}}, \langle l_{2}, Y \rangle_{\cal{H}})}(z_1, z_2)}{\partial z_1}$ and $\frac{\partial f_{(\langle l_{1}, X \rangle_{\cal{H}}, \langle l_{2}, Y \rangle_{\cal{H}})}(z_1, z_2)}{\partial z_2}$ exist. 

\noindent (A3) The random elements $X$ and $Y$ are norm bounded, i.e., $||X||_{\cal{H}}$ and $||Y||_{\cal{H}}$ are bounded random variables. Moreover, the joint probability density function of $\langle l_{1}, X \rangle_{\cal{H}}$ and $\langle l_{2}, Y \rangle_{\cal{H}}$ is uniformly bounded.

\noindent (A4) The first and the second derivatives of the probability density functions of the random variables $\langle l_{1}, X \rangle_{\cal{H}}$ and $\langle l_{2}, Y \rangle_{\cal{H}}$ are uniformly bounded. 

\noindent (A5) For $i = 1$ and 2, the kernels $k_i : \mathbb{R}\rightarrow A\subset\mathbb{R}^{+}$ are such that $\int k_{i}(m) dm = 1$, $\int m k_{i} (m) dm < \infty$, $\int k_{i}^{2} (m) dm < \infty$, $\int m^{2} k_{i} (m) dm < \infty$ and $\int m \{k_{i}(m)\}^{2} dm < \infty$. Here $A$ is a bounded set. 

\noindent (A6) The bivariate kernel $k: \mathbb{R}^{2}\rightarrow A\subset\mathbb{R}^{+}$ is such that $\int m_{1} k(m_1, m_2) dm_1 dm_2 < \infty$, $\int m_{2} k(m_1, m_2) dm_1 dm_2 < \infty$,  $\int m_{1} \{k(m_1, m_2)\}^{2} dm_1 dm_2 < \infty$ and $\int m_{2} \{k(m_1, m_2)\}^{2} dm_1 dm_2 < \infty$. Here $A$ is a bounded set.

\begin{rmrk}
Condition (A1) indicates that for various choices of $h_{n}$, the asymptotic results described in Theorem \ref{asymptotic} holds, and such conditions are common across the literature in kernel density estimation (see, e.g., \cite{Silverman1998}). As an application of (A1), we have $n h_{n}^{2}\rightarrow\infty$ as $n\rightarrow\infty$; this observation shall be used repeatedly in our analysis. Condition (A2) is satisfied when the directional derivatives of the joint probability density function of $||X||_{\cal{H}}$ and $||Y||_{\cal{H}}$ exist, and  this is indeed a realistic assumption.  Further, when the infinite dimensional random elements $X$ and $Y$ are norm bounded, and the joint probability density function of the random variables $\langle l_{1}, X \rangle_{\cal{H}}$ and $\langle l_{2}, Y \rangle_{\cal{H}}$ is uniformly bounded (see Proposition \ref{propbounded}), then the marginal probability density functions of the random variables $\langle l_{1}, X \rangle_{\cal{H}}$ and $\langle l_{2}, Y \rangle_{\cal{H}}$ are also bounded. For instance, for well-known Gaussian process, the conditions (A3) and (A4) hold. Condition (A5) implies certain mild restriction on the tail behaviour of the chosen kernels. However, for compact support of the kernel, there is no need to put restrictions on the tail behaviour of the chosen kernels as the integrations involving kernels are finite as long as the the kernels are continuous.  Condition (A6) asserted some integrability assumptions on the chosen bivariate kernel. In this case also, if the support of the chosen bivariate kernel is compact, the continuity of the chosen bivariate kernel is good enough to satisfy the integrability assumptions in condition (A6).  
\end{rmrk}

\begin{theorem}\label{asymptotic}
Let us denote \begin{equation}\label{practicestatistic}
\begin{aligned}
T_{n}^{G, L} &= 
\sup_{\substack{\theta_{j, i}\in (-\frac{\pi}{2}, \frac{\pi}{2})\\ j = 1, 2; i = 1, \ldots, M - 1\\ \theta_{j, M - 2}\in (-\pi, \pi)\\ j = 1, 2\\ s\in\{s_1, \ldots, s_L\}\\ t\in\{t_1, \ldots, t_L\}}}
\left|\frac{1}{n h_{n}^{2}}\sum\limits_{i = 1}^{n}k\left(\frac{\langle \hat{l}_{1}^{K}, X_{i}\rangle_{\cal{H}} - s}{h_{n}}, \frac{\langle \hat{l}_{2}^{K}, Y_{i}\rangle_{\cal{H}} - t}{h_{n}}\right)\right.\\ 
&\left. - \frac{1}{n^2 h_{n}^{2}}\sum\limits_{i, j = 1}^{n}k_{1}\left(\frac{\langle \hat{l}_{1}^{K}, X_{i}\rangle_{\cal{H}} - s}{h_{n}}\right)k_{2}\left(\frac{\langle \hat{l}_{2}^{K}, Y_{j}\rangle_{\cal{H}} - t}{h_{n}}\right)\right|,  
\end{aligned}
\end{equation} where $s_{i} = -G + i \frac{2G}{L}$ and $t_{i} = -G + i \frac{2G}{L}$ for $i = 1, \ldots, L$ are equi-distributied points in $[-G, G]$, and $T_{n}$ is the same as defined in \eqref{statistic}. Then, under (A1), (A5) and (A6), $$\sqrt{n}h_{n}\left\{\displaystyle\lim_{G\rightarrow\infty}\displaystyle\lim_{L\rightarrow\infty} (T_{n} - T_{n}^{G, L})\right\}\stackrel{p}\rightarrow 0$$ as $n\rightarrow\infty$. 

Moreover, under conditions (A1)--(A6) and when $T = \frac{\lambda}{\sqrt{n} h_{n}}$ for some $\lambda > 0$, $\sqrt{n} h_{n} T_{n}^{G, L} - \lambda$ converges weakly to the distribution of $\sup\limits_{s, t \in \{-G + i \frac{2G}{L}\mid i = 1, \ldots, L\}}\{|Z_{s, t}|\}$ as $n\rightarrow\infty$, $K\rightarrow\infty$ and $M\rightarrow\infty$ (see \eqref{finall} for the definitions of $K$ and $M$), where $Z_{s, t}$ is a random variable associated with normal distribution with mean\\ $= \sqrt{c}\left\{\sum\limits_{i = 1}^{2}\left(\frac{\partial f_{Z_{1}, Z_{2}}(s, t)}{\partial Z_{i}}\right)\int m_{i} k(m_1, m_2)dm_1 dm_2\right.$\\
$\left. - f_{Z_{1}}(s)f^{'}_{Z_2}(t)\int mk_{2}(m) dm - f_{Z_{2}}(t)f^{'}_{Z_1}(s)\int mk_{1}(m) dm\right\}$ \\and variance $= f_{Z_1, Z_2}(s, t)\int k^{2}(m_1, m_2) dm_1 dm_2$. Here $c = \displaystyle\lim_{n\rightarrow\infty} nh_{n}^{4}$ and $(Z_1, Z_2) = (\langle l_1, X\rangle_{\cal{H}}, \langle l_2, Y\rangle_{\cal{H}})$.
\end{theorem}

\begin{rmrk}
The first part in the statement of Theorem \ref{asymptotic} indicates that $T_{n}^{G, L}$ is an appropriate approximation of $T_{n}$ for sufficiently large $G$ and $L$ after normalized by $\sqrt{n} h_{n}$. Here $G$ controls the length of the interval of the time parameters, and $L$ controls the number of time points after discretization. The motivation of constructing $T_{n}^{G, L}$ will be discussed in detail in Section \ref{CT}. Moreover, with the same normalization factor $\sqrt{n} h_{n}$, $\sqrt{n} h_{n}T_{n}^{G, L} - \lambda$ converges weakly to a certain random variable associated with the functional of multivariate normal distribution. All together, one can approximately compute the asymptotic power of the test based on $T_{n}$ using the assertion in Theorem \ref{asymptotic} 
\end{rmrk}

\begin{rmrk}
Observe that second part of Theorem \ref{asymptotic} holds when $M\rightarrow\infty$ and $K\rightarrow\infty$, and it implies from \eqref{finall} that the stated asymptotic result is valid for uniform choices of $l_1$ and $l_2$. In the same spirit, $T_{n}^{G, L}$ approximates $T_{n}$ when $G\rightarrow\infty$ and $L\rightarrow\infty$ as we discussed in the previous remark, and it implies that the result is valid for uniform choices of time parameters $s$ and $t$. In all, the infinite limits of $M$, $K$, $G$ and $L$ connect the infinite dimensional space and its finite discretization.  
\end{rmrk}

\begin{rmrk}
Observe that Theorem \ref{asymptotic} is stated when $T = \frac{\lambda}{\sqrt{n} h_{n}}$, and since $\sqrt{n}h_{n}\rightarrow\infty$ as $n\rightarrow\infty$ (see Condition (A1)), it indicates that the statement of this theorem is valid for local alternatives. As a result, the assertion of this theorem enables us to compute the approximated asymptotic power of the test based on $T_{n}$. Further, as a special case, when $\lambda = 0$, the asymptotic distribution of $T_{n}$ under null hypothesis (i.e., $T = 0$) follows from Theorem \ref{asymptotic}, and it is formally stated in Corollary \ref{asymptoticnull}.    
\end{rmrk}

\begin{crlre}\label{asymptoticnull}
Under conditions (A1)-(A6) and when $T = 0$ (i.e., under null hypothesis), $\sqrt{n} h_{n}T_{n}^{G, L}$ converges weakly to the distribution of $\sup\limits_{s, t \in \{-G + i \frac{2G}{L}\mid i = 1, \ldots, L\}}\{|Z_{s, t}|\}$ as $n\rightarrow\infty$, $K\rightarrow\infty$ and $M\rightarrow\infty$, where $Z_{s, t}$ denotes a random variable associated with normal distribution with mean and variance, which are the same as described in the statement of Theorem \ref{asymptotic}.
\end{crlre}

Corollary \ref{asymptoticnull} asserts the asymptotic distribution of $T_{n}^{G, L}$ under the null hypothesis, and therefore, one can estimate the critical value of the test based on $T_{n}$ using the asymptotic distribution of $T_{n}^{G, L}$ stated in Corollary \ref{asymptoticnull} as $T_{n}$ and $T_{n}^{G, L}$ can be made arbitrary close to each other for sufficiently large $G$ and $L$ (see the first part of Theorem \ref{asymptotic}). All in all, one can explicitly express the power of the test based on $T_{n}^{G, L}$ when $T = t_{0}$ using Corollary \ref{asymptoticnull} and Theorem \ref{asymptotic}. In notation, if $\hat{c}_{\alpha}$ is such that $P_{H_{0}}[\sqrt{n} h_{n}T_{n}^{G, L} > \hat{c}_{\alpha}]\rightarrow \alpha$ as $n\rightarrow\infty$, then the power of the test based on $T_{n}^{G, L}$ under $H_{1}$ (i.e., when $T = t_{0}$) is given by 
\begin{equation}\label{power}
P_{H_{1}}\left[\sqrt{n} h_{n}T_{n}^{G, L} > \hat{c}_{\alpha}\right] = P\left[\sqrt{n}h_{n}(T_{n}^{G, L} - t_{0})>\hat{c}_{\alpha} - {\sqrt{n}h_{n}}t_{0}\right],
\end{equation} and the asymptotic power of the test based on $T_{n}^{G, L}$ under local alternatives (e.g., $T = {\frac{\lambda}{\sqrt{n} h_{n}}}$ as mentioned in Theorem \ref{asymptotic}) is given by 
\begin{equation}\label{asymptoticpower}
P_{\left(T = {\frac{\lambda}{\sqrt{n}h_{n}}}\right)}
\left[\sqrt{n h_{n}}T_{n}^{G, L} > \hat{c}_{\alpha}\right] = P\left[\sqrt{nh_{n}}T_{n}^{G, L} - \lambda>\hat{c}_{\alpha} - \lambda\right].
\end{equation} Moreover, using \eqref{power}, one can establish the consistency of the test based on $T_{n}^{G, L}$, which is stated in Proposition \ref{consistency}. 
\begin{proposition}\label{consistency}
Let $\hat{c}_{\alpha}$ is such that $P_{H_{0}}[{\sqrt{n}h_{n}}T_{n}^{G, L} > \hat{c}_{\alpha}]\rightarrow\alpha$ as $n\rightarrow\infty$, where $\alpha\in (0, 1)$ is a preassigned constant. Then, the test based on $T_{n}^{G, L}$, which rejects $H_{0}$ if ${\sqrt{n} h_{n}}T_{n}^{G, L}\geq\hat{c}_{\alpha}$, is a consistent test. In other words, $P_{H_{1}}[\sqrt{n} h_{n}T_{n}^{G, L}\geq\hat{c}_{\alpha}]\rightarrow 1$ as $n\rightarrow\infty$. 
\end{proposition} 
Proposition \ref{consistency} asserts that the test based on $T_{n}$ is expected to perform well in terms of power for a sufficiently large sample, as $T_n$ and $T_n^{G, L}$ can be made arbitrarily close for large $G$ and $L$ (see Proposition \ref{GL}). Next, we study the performance of the test based on $T_{n}$ through the asymptotic power study.  

\subsubsection{Asymptotic Power Study}\label{APS}
Recall that Theorem \ref{asymptotic} states the asymptotic distribution of $T_{n}^{G, L}$ under local alternatives, i.e., when $T = {\frac{\lambda}{\sqrt{n} h_n}}$ ($\lambda\geq 0$), and one can compute the asymptotic power of the test based on $T_{n}$ for various choices of the random processes in separable Hilbert space and $\lambda$. In the course of study, the critical value with $\alpha\%$ level of significance, i.e., $\hat{c}_{\alpha}$ is computed using the result described in Corollary \ref{asymptoticnull}, and using the obtained $\hat{c}_{\alpha}$, one can approximate the asymptotic power of the test based on $T_{n}$ using \eqref{asymptoticpower}. 

We consider following three examples : 

\noindent {\bf Example 1} : Let $(X, Y)$ be $L_{2}([0, 1])\times L_{2}([0, 1])$ valued random element. Suppose that $X(t)$ ($t\in [0, 1]$) has the same distribution as $U\exp(t)$, where $U$ follows uniform distribution over $[0, 1]$, and $Y (t) = ||X||_{L_{2}([0, 1])}\times B(t)$, where $B(t)$ ($t\in [0, 1]$) is a Gaussian process with  $E(B(t))= 0$ and $E(B(t)B(s)) = \min(t, s)$ for all $t\in [0, 1]$ and $s\in [0, 1]$. 

\noindent {\bf Example 2} : Let $(X, Y)$ be $L_{2}([0, 1])\times L_{2}([0, 1])$ valued random element. Suppose that $X(t)$ ($t\in [0, 1]$) has the same distribution as $U\exp(t)$, where $U$ follows uniform distribution over $[0, 1]$, and $Y (t) = ||X||_{L_{2}([0, 1])}\times B_{1}(t)$, where $B_{1}(t)$ ($t\in [0, 1]$) is fractional Brownian motion with Hurst index $= .25$, i.e., $E(B_{1}(t)) = 0$ and $E(B_{1}(t)B_{1}(s)) = \frac{1}{2}(|t|^{0.5} + |s|^{0.5} - |t - s|^{0.5})$ for all $t\in [0, 1]$ and $s\in [0, 1]$. 

\noindent {\bf Example 3} : Let $(X, Y)$ be $L_{2}([0, 1])\times L_{2}([0, 1])$ valued random element. Suppose that $X(t)$ ($t\in [0, 1]$) has the same distribution as $U\exp(t)$, where $U$ follows uniform distribution over $[0, 1]$, and $Y (t) = ||X||_{L_{2}([0, 1])}\times B_{2}(t)$, where $B_{2}(t)$ ($t\in [0, 1]$) is fractional Brownian motion with Hurst index $= .75$, i.e., $E(B_{2}(t)) = 0$ and $E(B_{2}(t) B_{2}(s)) = \frac{1}{2}(|t|^{1.5} + |s|^{1.5} - |t - s|^{1.5})$ for all $t\in [0, 1]$ and $s\in [0, 1]$.

\begin{figure}
	\includegraphics[width= 2in, height = 2.5in]{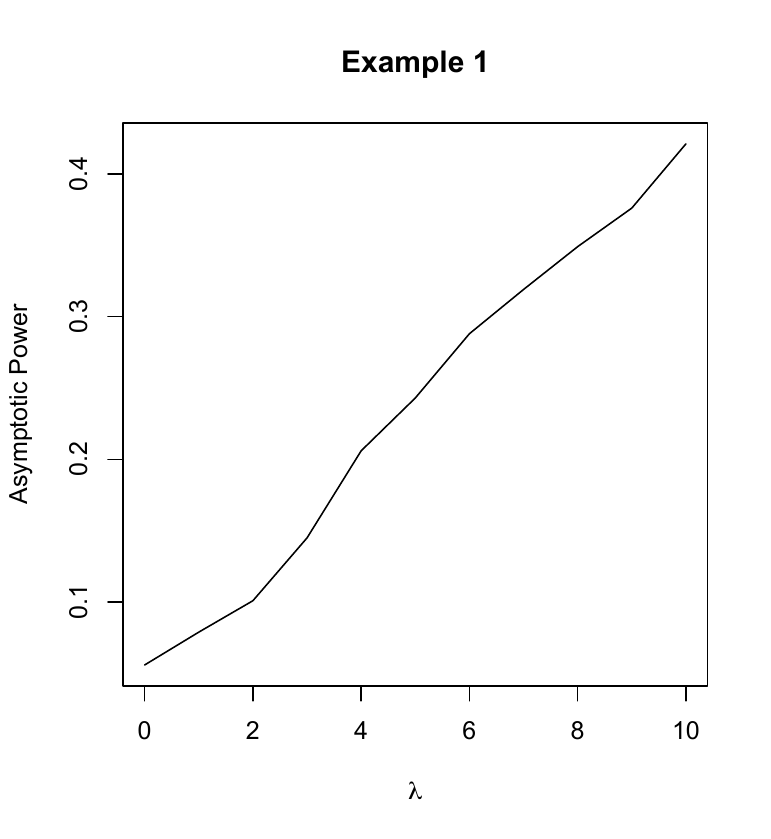}
	\includegraphics[width= 2in, height = 2.5in]{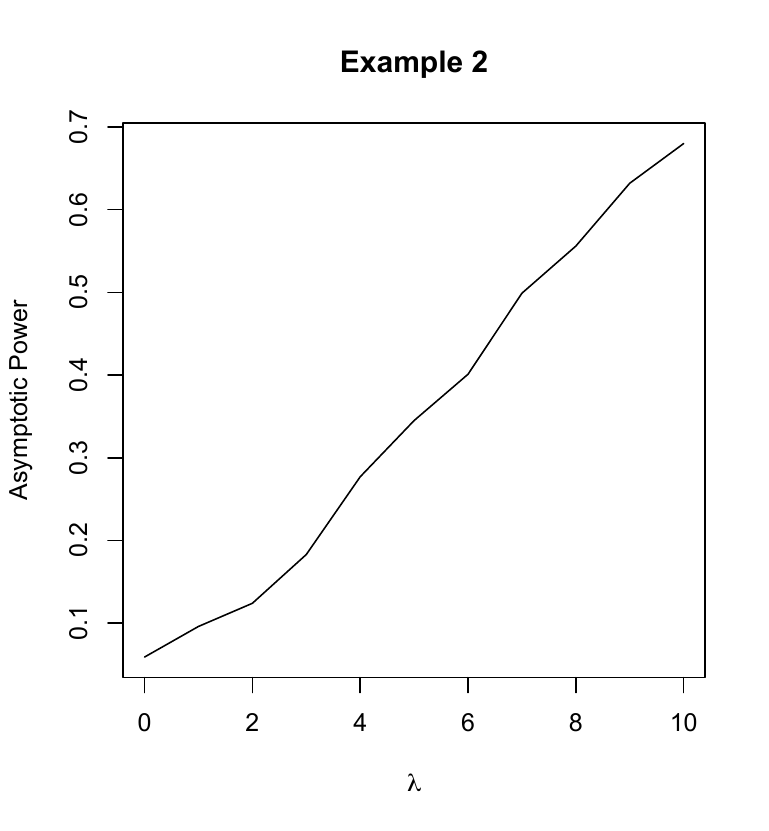}\includegraphics[width= 2in, height = 2.5in]{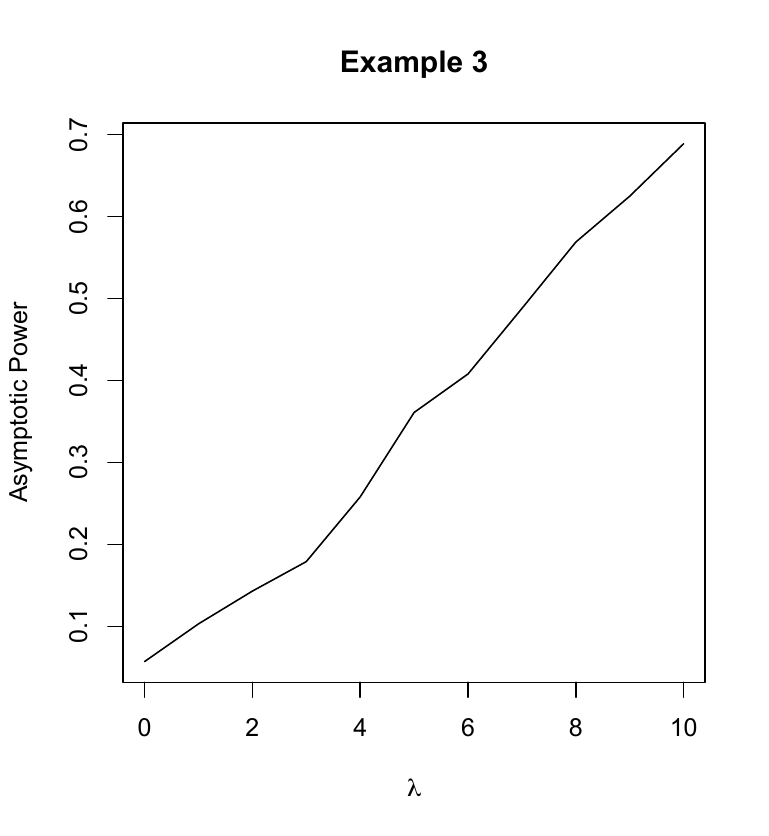}
	\caption{\it Asymptotic power of the proposed test for Examples 1, 2 and 3 at 5\% level of significance when $\lambda$ varies from 0 to 10. 
		} \label{Fig_AsP}
\end{figure}

Note that in order to compute the asymptotic power of the test, one first needs to compute the critical value, which can be estimated using the result stated in Corollary \ref{asymptoticnull}, i.e., when $T = 0$, and it is equivalent to the fact that $X$ and $Y$ are independent random elements in ${\cal{H}}$ (see Theorem \ref{characterization}). Now, in order to satisfy this statistical independence between $X$ and $Y$, we first consider $X(t)$ and $Y(t)$ are two independent random elements in $L_{2}([0, 1])$ associated with a standard Brownian motion. Afterwards, the terms involving the joint and the marginal probability density functions (i.e.,  $f_{(\langle l_{1}, X \rangle_{{\cal{H}}}, \langle l_{2}, Y \rangle_{{\cal{H}}})}(., .), f_{\langle l_{1}, X \rangle_{{\cal{H}}}} (.), ~\mbox{and}~ f_{\langle l_{2}, Y \rangle_{{\cal{H}}}} (.)$) appeared in the mean and the variance of the asymptotic distribution of $T_{n}$ under $H_{0}$ (see Corollary \ref{asymptoticnull}) are estimated by bivariate and marginal kernel density function. For the bivariate kernel density function, we consider the two fold product of the marginal kernel density function. In this study, for estimating the marginal probability density function, we use the well-known Epanechnikov Kernel because of its optimal properties (see, e.g., \cite{Silverman1998}). Finally, $\alpha\%$ (denoted as $\hat{c}_{\alpha}$) critical value is computed as the $(1 - \alpha)$-th quantile of the normal distribution described in Corollary \ref{asymptoticnull}. Precisely speaking,
\begin{equation}\label{estcritical}
\hat{c}_{\alpha} = G_{H_{0}}^{-1}(\alpha),    
\end{equation} where $G_{H_{0}}$ is the cumulative distribution function of ${\sqrt{n} h_{n}}T_{n}^{G, L}$ (see Corollary \ref{asymptoticnull} for the explicit form of $G_{H_{0}}$). In the study, we consider $G = 20$ and $L = 10$.

Now, for Examples 1, 2 and 3, the approximated asymptotic powers of the test based on $T_{n}$ are illustrated in the diagrams in Figure \ref{Fig_AsP} for various choices of $\lambda$. In general, it is observed from all three diagrams that the asymptotic power of the test based on $T_{n}$ increases as $\lambda$ increases. The increasing feature of asymptotic power with respect to $\lambda$ for all three examples has been seen in view of the fact that $X$ and $Y$ are {\it not} independent random elements and as the value of $T$ is deviating more from the null hypothesis as $\lambda$ is deviating from zero. Moreover, among Examples 1, 2 and 3, observe that the asymptotic powers for Examples 2 and 3 are more than that for Example 1 since in Examples 2 and 3, $Y(t)$ at various time points $t\in [0, 1]$ are correlated whereas in Example 1, $Y(t)$ at various time points $t\in [0, 1]$ are uncorrelated. 

\section{Finite Sample Level and Power Study}\label{FSLPS}
The asymptotic power study in Section \ref{APS} indicates that the power of the test based on $T_{n}$ is fairly good when the sample size is large enough. Besides, Proposition \ref{consistency} asserts the consistency of the test, which also indicates that the proposed test will perform well in terms of power for a sufficiently large sample. Overall, these facts only guarantee the expected good performance of the proposed test when the sample size is large enough. Here, we now want to see the performance of the proposed test when the sample size is small or moderately small/large, and in order to carry out this study, one needs to compute $T_{n}$ for a given data, which is described in the following. 

\subsection{Computation of $T_{n}$}\label{CT}
Observe that the exact distribution of $T_{n}$ is intractable because of its complex structure, and for that reason, one cannot compute the critical value from the inverse of the exact distribution of $T_{n}$. In order to overcome this issue, one requires to approximate the sampling distribution of $T_{n}$ using Monte-Carlo simulation, and in this procedure, the computation of $T_{n}$ for a given sample is requisite. Further, for the same reason provided above, one needs to compute $T_{n}$ for a given sample for computing the power. 

Recall $T_{n}$ from \eqref{statistic} : 
\begin{eqnarray*}
T_{n} &= 
\displaystyle\sup_{\substack{\theta_{j, i}\in (-\frac{\pi}{2}, \frac{\pi}{2})\\ j = 1, 2; i = 1, \ldots, M - 2\\ \theta_{j, M - 1}\in (-\pi, \pi)\\ j = 1, 2\\ s\in\mathbb{R}, t\in\mathbb{R}}}
\left|\frac{1}{n h_{n}^{2}}\sum\limits_{i = 1}^{n}k\left(\frac{\langle \hat{l}_{1}^{K}, X_{i}\rangle_{\cal{H}} - s}{h_{n}}, \frac{\langle \hat{l}_{2}^{K}, Y_{i}\rangle_{\cal{H}} - t}{h_{n}}\right)\right.\\
&\left.- \frac{1}{n^2 h_{n}^2}\sum\limits_{i, j = 1}^{n}k_{1}\left(\frac{\langle \hat{l}_{1}^{K}, X_{i}\rangle_{\cal{H}} - s}{h_{n}}\right)k_{2}\left(\frac{\langle \hat{l}_{2}^{K}, Y_{j}\rangle_{\cal{H}} - t}{h_{n}}\right)\right|.
\end{eqnarray*}
Note that in the expression in $T_{n}$, the supremum is taken over finitely many choices of $\theta_{j, i}$ ($j = 1$ and 2, $i = 1, \ldots, M - 1$), which are involved in $\hat{l}_{1}^{K}$ and $\hat{l}_{2}^{K}$ and over $s\in\mathbb{R}$ and $t\in\mathbb{R}$. As $\mathbb{R}$ in an interval with infinite length, the computation of supremum on $s$ and $t$ over $\mathbb{R}$ becomes unmanageable, and in order to overcome this issue, we take a sufficiently large interval $[-G, G]$ and do the grid search on $\{s_1, \ldots, s_L\}$ and $\{t_1, \ldots, t_L\}$, where $s_{i}\in [-G, G]$ and $t_{i}\in [-G, G]$ ($i = 1, \ldots, L$) are equally spaced points. This approximation procedure makes the computation doable in practice since domain of both variables $s$ and $t$ become finite. Therefore, in practice, we will compute the following.  

\begin{equation}\label{practicestatistic}
\begin{aligned}
T_{n}^{G, L} &= 
\sup_{\substack{\theta_{j, i}\in (-\frac{\pi}{2}, \frac{\pi}{2})\\ j = 1, 2; i = 1, \ldots, M - 1\\ \theta_{j, M - 2}\in (-\pi, \pi)\\ j = 1, 2\\ s\in\{s_1, \ldots, s_L\}\\ t\in\{t_1, \ldots, t_L\}}}
\left|\frac{1}{n h_{n}^{2}}\sum\limits_{i = 1}^{n}k\left(\frac{\langle \hat{l}_{1}^{K}, X_{i}\rangle_{\cal{H}} - s}{h_{n}}, \frac{\langle \hat{l}_{2}^{K}, Y_{i}\rangle_{\cal{H}} - t}{h_{n}}\right)\right.\\ 
&\left. - \frac{1}{n^2 h_{n}^{2}}\sum\limits_{i, j = 1}^{n}k_{1}\left(\frac{\langle \hat{l}_{1}^{K}, X_{i}\rangle_{\cal{H}} - s}{h_{n}}\right)k_{2}\left(\frac{\langle \hat{l}_{2}^{K}, Y_{j}\rangle_{\cal{H}} - t}{h_{n}}\right)\right|,  
\end{aligned}
\end{equation} where $G\in\mathbb{R}$ and $L\in\mathbb{Z}^{+}$ are sufficiently large. The proposition \ref{GL} along with the first part of Theorem \ref{asymptotic} guide our choice of $T_{n}^{G, L}$ as an approximation of $T_{n}$. 

\begin{ppsn}\label{GL}
For all $n\geq 1$, under (A5) and (A6), $\displaystyle\lim_{G \to \infty} \displaystyle\lim_{L \to \infty} (T_{n} - T_{n}^{G, L})\stackrel{a.s.} =  0$, where $T_{n}$ is defined in \eqref{statistic}, and $T_{n}^{G, L}$ is defined in \eqref{practicestatistic}.
\end{ppsn} 


\subsubsection{Estimated Level and Power}\label{ELP}
Here, we discuss how one can estimate the level and the power of the test based on $T_{n}$. In this study, we choose $M = 10$ and $K = 20$ as we have observed that any larger values than the aforementioned chosen values does not much change the result. For details about the variables $M$ and $K$, recall Proposition \ref{lk} and the  discussion before the proposition. Besides, as Proposition \ref{GL} indicates, the choices of $G$ and $L$ are also an issue of concern. We observe that $G = 20$ and $L =10$ are reasonably good choices in terms of computational complexity and the precision of this finite approximation. Throughout the study, we consider $k_1$ and $k_2$ in the expression in $T_{n}$ (see \eqref{statistic}) as Epanechnikov kernel (see \cite{Silverman1998}) for its optimal property, and the bivariate kernel $k$ is considered as the two folds product of Epanechnikov kernels. As the bandwidth, we choose $h_{n} = c n^{-\frac{1}{6}}$, where $c$ is estimated by cross validation technique (see, e.g., \cite{Duong2005} and \cite{Hall1992}). In the course of this study, the experiments are carried out for $n = 20$, $50$ and $100$. 

In order to carry out the critical value of the test, we generate two independent data sets ${\cal{X}} = (X_{1}, \ldots, X_{n})$ and ${\cal{Y}} = (Y_{1}, \ldots, Y_{n})$ from standard Brownian motion, and replicate this experiment $1000$ times. Note that as ${\cal{X}}$ and ${\cal{Y}}$ are independent, the null hypothesis $H_{0} : X\perp Y$ ($\Leftrightarrow T = 0$) is satisfied. Now, let $T_{n, i}$ be the value of $T_{n}$ for the $i$-th replicate $(i = 1, \ldots, 1000)$, and these 1000 values of $T_{n, i}$ provide the estimated distribution of $T_{n}$. Therefore, $(1 - \alpha)$-th ($\alpha\in (0, 1)$) quantile of $T_{n, i}$ $(i = 1, \ldots, 1000)$, denoted as $\hat{c}_{\alpha}$ can be considered as the estimated critical value. Now, in order to estimate the level of significance, we generate two {\it independent} data sets, namely, ${\cal{X}}^{'} = (X_{1}^{'}, \ldots, X_{n}^{'})$ and ${\cal{Y}}^{'} = (Y_{1}^{'}, \ldots, Y_{n}^{'})$ and replicate this experiment $M$ times. Suppose that $T_{n, i}$ is the value of $T_{n}$ for the $i$-th replication, and the estimated level is defined as $\frac{1}{M}\sum\limits_{i = 1}^{M} 1_{(T_{n, i} > \hat{c}_{\alpha})}$. Next, in order to estimate the power of the test, we generate two  {\it dependent} data sets, namely, ${\cal{X}}^{''} = (X_{1}^{''}, \ldots, X_{n}^{''})$ and ${\cal{Y}}^{''} = (Y_{1}^{''}, \ldots, Y_{n}^{''})$ and replicate this experiment $M^{'}$ times. Suppose that $T_{n, j}$ is the value of $T_{n}$ for the $j$-th replication, and the estimated power is defined as $\frac{1}{M^{'}}\sum\limits_{j = 1}^{M^{'}} 1_{(T_{n, j} > \hat{c}_{\alpha})}$. 

To estimate the level, we consider the following three examples. 

\noindent {\bf Example 4 :} Two {\it independent} data sets ${\cal{X}} = \{X_{1}, \ldots, X_{n}\}$ and ${\cal{Y}} = \{Y_1, \ldots, Y_{n}\}$ are generated from standard Brownian motion. 

\noindent {\bf Example 5 :} Two {\it independent} data sets ${\cal{X}} = \{X_{1}, \ldots, X_{n}\}$ and ${\cal{Y}} = \{Y_1, \ldots, Y_{n}\}$ are generated from fractional Brownian motion with Hurst index $= 0.25$. 
 
\noindent {\bf Example 6 :} Two {\it independent} data sets ${\cal{X}} = \{X_{1}, \ldots, X_{n}\}$ and ${\cal{Y}} = \{Y_1, \ldots, Y_{n}\}$ are generated from fractional Brownian motion with Hurst index $= 0.75$. 

\begin{table}[h!]
	
	\begin{center}		
		
		\begin{tabular}{ccccc}\hline
			model & $n = 20$ & $n = 50$ & $n = 100$ & $n = 500$\\ \hline 
			Example 4 ($\alpha = 5\%$) & {$0.062$} & {$0.060$} & {$0.057$} & {$0.051$}\\ \hline
			Example 4 ($\alpha = 10\%$) & {$0.121$} & {$0.120$} & {$0.114$} & {$0.107$}\\ \hline
		    Example 5 ($\alpha = 5\%$) & {$0.063$} & {$0.058$} & {$0.054$} & {$0.052$}\\ \hline
		    Example 5 ($\alpha = 10\%$) & {$0.122$} & {$0.121$} & {$0.110$} & {$0.106$}\\ \hline
		    Example 6 ($\alpha = 5\%$) & {$0.064$} & {$0.062$} & {$0.060$} & {$0.056$}\\ \hline
		    Example 6 ($\alpha = 10\%$) & {$0.125$} & {$0.124$} & {$0.119$} & {$0.111$}\\ \hline
		\end{tabular}
	\end{center}
	\caption{\it The estimated size of the proposed test for different sample sizes $n$. The levels of significance, i.e., $\alpha$ are $5\%$ and $10\%$.}
	
	\label{tab1}
	
\end{table}

\noindent In order to generate the data from the fractional/standard Brownian motion, we generate the data corresponding to a certain large dimensional multivariate normal distribution. The estimated levels for Examples 4, 5 and 6 are reported in Table \ref{tab1} for the sample sizes $n = 20$, 50, 100 and 500 at 5\% and 10\% level of significance. The reported values indicate that the estimated level never deviated more than 1.5\% from desired level of significance (i.e., 5\% and 10\%) when the sample sizes are 20 and 50 and not deviated more than 1\% when the sample sizes are 100 and 500. 

\begin{table}[h!]
	
	\begin{center}		
		
		\begin{tabular}{ccccc}\hline
			model & $n = 20$ & $n = 50$ & $n = 100$ & $n = 500$\\ \hline 
			Example 1 ($\alpha = 5\%$) & {$0.342$} & {$0.506$} & {$0.678$} & {$0.791$}\\ \hline
			Example 1 ($\alpha = 10\%$) & {$0.399$} & {$0.601$} & {$0.800$} & {$0.882$}\\ \hline
			Example 2 ($\alpha = 5\%$) & {$0.445$} & {$0.624$} & {$0.785$} & {$0.911$}\\ \hline
		    Example 2 ($\alpha = 10\%$) & {$0.503$} & {$0.632$} & {$0.867$} & {$0.949$}\\ \hline
		    Example 3 ($\alpha = 5\%$) & {$0.427$} & {$0.641$} & {$0.805$} & {$0.902$}\\ \hline
		    Example 3 ($\alpha = 10\%$) & {$0.517$} & {$0.666$} & {$0.841$} & {$0.956$}\\ \hline
		
		\end{tabular}
	\end{center}
	\caption{\it The estimated power of the proposed test for different sample sizes $n$. The level of significance, i.e., $\alpha$ are $5\%$ and 10\%.}
	
	\label{tab2}
	
\end{table}

To estimate the power, we generate the data $(X_{1}, \ldots, X_{n})$ and $(Y_1, \ldots, Y_n)$ from the distributions as described in Examples 1, 2 and 3 in Section \ref{APS}.  At 5\% and 10\% level of significance, the estimated powers for those examples are reported in Table \ref{tab2} when the sample sizes are $20$, 50 and 100. In all these cases, we observe that the power increases as the sample size increases. Besides, as also observed in the asymptotic power study (see Section \ref{APS}), the proposed test becomes more powerful when the data is generated from the distribution described in Examples 2 and 3 than that of Example 1, and this fact is expected since the correlation structure involved in the fractional Brownian motion having Hurst parameter $\neq 0.5$.

A few more simulation studies are carried out in Appendix B.

\section{Real Data Analysis}\label{RDA}
\begin{figure}
	\includegraphics[width= 3.5in, height = 3.5in]{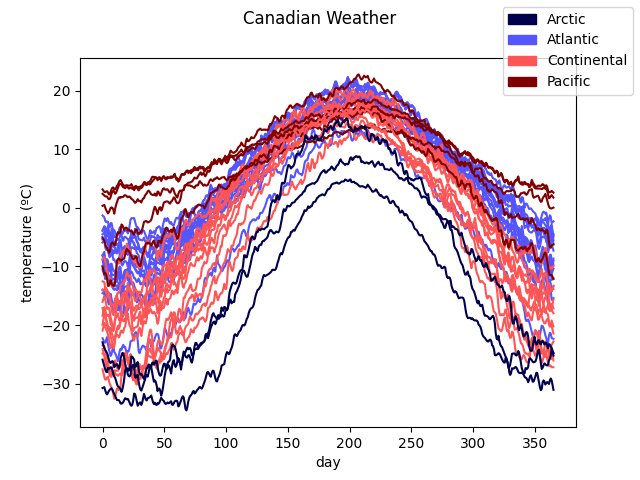}
	\includegraphics[width= 3.5in, height = 3.5in]{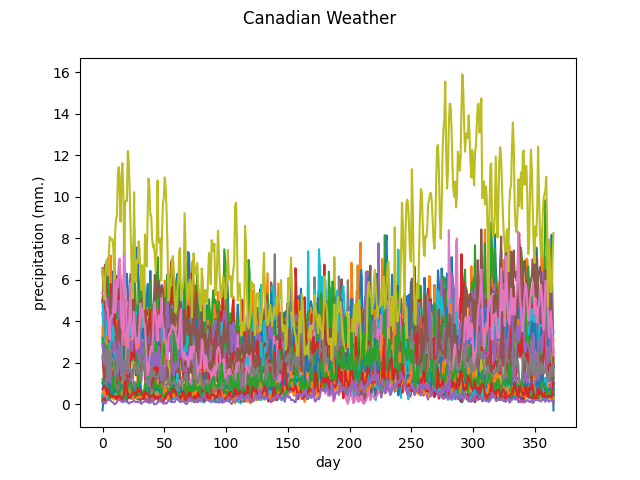}
	\caption{\it  
		\label{CanadianWeather}
		The left diagram  plots the average daily temperature and the right diagram plots the average daily precipitation. Source : \url{https://fda.readthedocs.io/en/latest/auto_examples}} 
\end{figure}

Here, we implement our proposed test on a real data, which consists of daily temperature and precipitation at $n = 35$ locations in Canada averaged over 1960 to 1994. \cite{Ramsey2005} studied this data set in the context of functional principal component analysis, and the soft version of the data may be found in \url{https://climate.weather.gc.ca/historical_data/search_historic_data_e.html}. The dimension of the data equals 365, which is much larger than the sample size $= 35$ of the data, and therefore, embedding such a high dimensional data into a specific Hilbert space is legitimate enough. For this data set, suppose that $X_{i}$ ($i = 1, \ldots, 35$) is the average daily temperature for each day of the year at the $i$-th location, and $Y_{i}$ is the base 10 logarithm of the corresponding average precipitation. Here we consider $X$ and $Y$ as $L_{2}([0, 1])$-valued random elements, where $365$ days are considered as the equally spaced 365 time points over $[0, 1]$. The curves associated with $X$ and $Y$ are demonstrated in the diagrams in Figure \ref{CanadianWeather}.

In order to check whether $X$ and $Y$ are independent random elements or not, we compute the $p$-value of the test based on $T_{n}$ using Bootstrap resampling procedure. At the beginning, to compute $T_{n}$ (precisely speaking, $T_{n}^{G, L}$), we choose tuning parameters $M = 10$, $K = 20$, $G = 20$ and $L = 10$ as this issue has been discussed at the beginning of Section \ref{ELP}. Similarly, we choose the bandwidth $h_{n} = c n^{-\frac{1}{6}}$, where $c$ is estimated by the well-known cross validation technique (see, e.g., \cite{Duong2005}), and as said before, $k_1$ and $k_2$ are chosen as Epanechnikov kernel, and the bivariate kernel $k$ is considered as the two fold product of Epanechnikov kernels. Based on these aforesaid choices, we generate 500 Bootstrap resamples with size $= 35$ from the original data $(X_1, Y_1), \ldots, (X_{35}, Y_{35})$, and suppose that $t_{0}$ is the value of $T_{n}$ for the original data $(X_1, Y_1), \ldots, (X_{35}, Y_{35})$. Let $t_{j}$ $(j = 1, \ldots, 500)$ be the value of $T_{n}$ for the $j$-th resample, and the $p$-value is defined as $\frac{1}{500}\sum\limits_{j = 1}^{500} 1_{(t_{j} > t_{0})}$. This methodology gives us the $p$-value as $0.078$, which indicates the data does not favour the null hypothesis, i.e., the temperature curve and the precipitation curve in those 35 locations in Canada over the period from 1960 to 1994 are not statistically independent at 8\% level of significance. Even from the climate science point of view, it is expected that the temperature and the amount of precipitation are supposed to have some dependence structure, and hence, the result shown by the proposed test is reasonable enough. 

For this data, as the temperature and the precipitation are dependent according to our proposed test, one may be interested to know the estimated size and power of the test for this data set. In order to compute estimated size and power, we combine $(X_1, \ldots, X_{35})$ and $(Y_{1}, \ldots, Y_{35})$ and mix these 70 observations randomly. Let us denote the new sample as ${\cal{Z}} = (Z_1, \ldots, Z_{70})$. Now, partition ${\cal{Z}}$ into two parts, namely, ${\cal{Z}}_{1} = (Z_1, \ldots, Z_{35})$ and ${\cal{Z}}_{2} = (Z_{36}, \ldots, Z_{70})$, and note that ${\cal{Z}}_{1}$ and ${\cal{Z}}_{2}$ are two independent data sets as these are constructed by random mixing of the combined sample of ${\cal{X}}$ and ${\cal{Y}}$. Next, from ${\cal{Z}}_{1}$ and ${\cal{Z}}_{2}$, by Bootstrap methodology, we generate $M$ many resamples ${\cal{Z}}_{1, i}$ and ${\cal{Z}}_{2, i}$ ($i = 1, \ldots, M$) and compute the test statistic $T_{n}$ for each resample. Suppose that $t_{i}$ is the value of $T_{n}$ for the $i$-th resample, and then, $(1 - \alpha)$-th quantile of $(t_1, \ldots, t_M)$ is considered as the estimated critical value (denoted as $\hat{c}_{\alpha}$) at $\alpha\%$ level of significance. In order to estimate the size of the test, we generate $M_1$ many resamples ${\cal{Z}}_{1, i}$ and ${\cal{Z}}_{2, i}$ ($i = 1, \ldots, M_1$), and the estimated size can be defined as $\frac{1}{M_1}\sum\limits_{j = 1}^{M_1}1_{(t_{j} > \hat{c}_{\alpha})}$, where $t_{j}$ is the value of $T_{n}$ for $j$-th resample $(j = 1, \ldots, M_1)$. Next, in order to estimate the power, using Bootstrap methodology, we generate $M_2$ many resamples from the original data $(X_{1}, \ldots, X_{35})$ and $(Y_1, \ldots, Y_{35})$, and note that in each resample, $X$-data and $Y$-data are statistically dependent since original $(X_{1}, \ldots, X_{35})$ and $(Y_1, \ldots, Y_{35})$ are statistically dependent data sets. Suppose that $t_{k}$ ($k = 1, \ldots, M_2$) is the value of $T_{n}$ for the $k$-th resample, and the estimated power can be defined as $\frac{1}{M_2}\sum\limits_{k = 1}^{M_1}1_{(t_{k} > \hat{c}_{\alpha})}.$ In this study, we consider $M_1 = M_2 = 500$ and estimate the power and the size of the test when $\alpha = 0.05$. Using all these choices, we obtain the estimated size $= 0.058$ and the estimated power $= 0.723$. Overall, these facts indicate that the proposed test can achieve the nominal level of the test and poses good power when the random elements are statistically dependent. 

Two more real data are analysed using the proposed methodology in Appendix B. 

\section{Concluding Remarks}\label{CR}
This article studies the test for independence of two random elements taking values on an infinite dimensional space. In this test, the test statistic is formulated based on the sup norm (i.e., $L_{\infty}$ norm) distance between the joint probability density function of two dimensional certain projection of bivariate random element and product of marginal probability density functions of corresponding marginal random variables. The asymptotic distribution of the test statistic under certain local alternatives have been derived, and it has been observed that the proposed test performs well for various choices of local alternatives. Besides, the usefulness of the test is shown on well-known data sets, and simulation studies also indicate that the proposed test performs well under different scenarios. 

It follows from the proof of Theorem \ref{asymptotic} that the main crux of the proof is the derivation of the pointwise asymptotic properties of $$\frac{1}{nh_{n}^{{2}}}\sum\limits_{i = 1}^{n}k\left(\frac{\langle \hat{l}_{1}^{K}, X_{i}\rangle_{\cal{H}} - s}{h_{n}}, \frac{\langle \hat{l}_{2}^{K}, Y_{i}\rangle_{\cal{H}} - t}{h_{n}}\right) - \frac{1}{n^2 h_{n}^{2}}\sum\limits_{i, j = 1}^{n}k_{1}\left(\frac{\langle \hat{l}_{1}^{K}, X_{i}\rangle_{\cal{H}} - s}{h_{n}}\right)k_{2}\left(\frac{\langle \hat{l}_{2}^{K}, Y_{j}\rangle_{\cal{H}} - t}{h_{n}}\right)\,$$ and afterwards, the asymptotic distribution of $T_{n}^{G, L}$ follows from some arguments related to the continuous mapping theorem (see, e.g., \cite{Serfling1980}). In this context, we would like to point out that one may consider any other appropriate norm like $L_p$ norm, when $p\in [1, \infty)$. The study of  performance of the test based on $L_{p}$ norm may be of interest for future research.

Another issue may need to be discussed, which is related to proposed criterion $T$ (see \eqref{finalycriterion}). One may argue whether $T$ can be thought of as a measure of association or not. Note that under some technical conditions, one can argue that $T\leq M$, where $M > 0$ is a constant. Hence, it may be appropriate to consider $T^{'} = \frac{T}{M}$ as a measure of association. Observe that $T^{'}\in [0, 1]$, and for any two random elements $X$ and $Y$, $T^{'} (X, Y) = 0$ if and only if $X$ and $Y$ are independent random elements, which follows from the assertion in Theorem \ref{characterization}. In order to establish $T^{'}$ as a measure of association, one needs to characterize the case when $T^{'} (X, Y) = 1$, i.e., in other words, the readers may be interested to know the {\it perfect} dependence structure between the random elements $X$ and $Y$, which is not done in this work. The characterization of {\it perfect} dependence between $X$ and $Y$ when $T^{'} (X, Y) = 1$ may be of another interest for future research.

\noindent{\bf Acknowledgement:}
The first author acknowledges the research grant DST/INSPIRE/04/2017/002835, Government of India,  and the second author is supported by the research grant CRG/2022/001489, Government of India.

\section{Appendix A : Technical Details}\label{TD}

\noindent {\bf Proof of Lemma \ref{lemindep}:} Let $X$ and $Y$ be independent, and fix $l_1$ and $l_2 \in \cal{H}$. Observe that the mappings $\inpr{l_1}{\cdot}:\cal{H} \to \mathbb{R}$ and $\inpr{l_2}{\cdot}:\cal{H} \to \mathbb{R}$ defined by $x \mapsto \inpr{l_1}{x}$ and $x \mapsto \inpr{l_2}{x}$ are measurable.

\noindent Now, for all Borel subsets $A$ and $B$ of $\R$, we have 
\begin{align*}
    P(\inpr{l_1}{X} \in A, \inpr{l_2}{Y} \in B) &= P(X \in \inpr{l_1}{\cdot}^{-1}(A), Y \in \inpr{l_1}{\cdot}^{-1}(B))\\
    &= P(X \in \inpr{l_1}{\cdot}^{-1}(A)) P(Y \in \inpr{l_1}{\cdot}^{-1}(B))\\
    &~\mbox{(since $X$ and $Y$ are independent)}\\
    &= P(\inpr{l_1}{X} \in A) \p(\inpr{l_2}{Y} \in B).
\end{align*}
Therefore, the two random variables $\inpr{l_1}{X}$ and $\inpr{l_2}{Y}$ are independent.

Conversely, suppose that $\inpr{l_1}{X}$ and $\inpr{l_2}{Y}$ are independent for all $l_1$ and $l_2 \in \cal{H}$. Since, $\cal{H}$ is separable, consider an orthonormal basis $\{e_n: n = 1, 2, \cdots\}$. Observe that
\[X \overset{\cal{H}}{=} \sum_{n = 1}^\infty \inpr{X}{e_n} e_n, \quad Y \overset{\cal{H}}{=} \sum_{n = 1}^\infty \inpr{Y}{e_n} e_n.\] Here $A\stackrel{\cal{H}} = B$ denotes $||A - B||_{\cal{H}}= 0$ for any $A\in{\cal{H}}$ and $B\in{\cal{H}}$.
Since the two families of real valued random variables $\{\inpr{X}{e_n}: n = 1, 2, \cdots\}$ and $\{\inpr{Y}{e_n}: n = 1, 2, \cdots\}$ are independent, so are $X$ and $Y$. This completes the proof. \hfill\qed

\noindent The following lemma is useful in proving Proposition \ref{prilimarycriterionstatement}.

\begin{lemma}\label{independence-under-scalar-mult}
Let $(Z, W)$ denote an $\IR^2$ valued random vector, and fix non-zero real numbers $\alpha$ and $\beta$. Then, $(Z, W)$ is absolutely continuous if and only if $(\alpha Z, \beta W)$ is absolutely continuous. Moreover, $Z$ and $W$ are independent if and only if $\alpha Z$ and $\beta W$ are independent.
\end{lemma}

\noindent {\bf Proof of Lemma \ref{independence-under-scalar-mult}:} Consider the transformation $(x, y)^t \in \R^2 \mapsto (\alpha x, \beta y)^t \in \R^2$. Observe that the transformation is one-to-one and onto, and the Jacobian of the transformation is a non-zero constant. Hence, the first statement related to absolute continuity follows.

To prove the second part, let $A$ and $B$ be Borel subsets of $\R$. Then, $\frac{1}{\alpha}A$ (elementwise product) and $\frac{1}{\beta} B$ (elementwise product) are also Borel subsets of $\R$. Now,
\begin{align*}
    & \alpha Z \text{ and } \beta W \text{ are independent}\\
    &\iff \p(\alpha Z \in A, \beta W \in B) = \p(\alpha Z \in A)\ \p(\beta W \in B), \text{ for all Borel subsets } A, B\\
    &\iff \p(Z \in \tfrac{1}{\alpha} A, W \in \tfrac{1}{\beta} B) = \p(Z \in \tfrac{1}{\alpha} A)\ \p(W \in \tfrac{1}{\beta} B), \text{ for all Borel subsets } A, B\\
    &\iff \p(Z \in A^\prime, W \in B^\prime) = \p(Z \in A^\prime)\ \p(W \in B^\prime), \text{ for all Borel subsets } A^\prime, B^\prime\\
    &\iff Z \text{ and } W \text{ are independent}
\end{align*}
This completes the proof. \hfill \qed

\noindent {\bf Proof of Proposition \ref{prilimarycriterionstatement}:}
Let $X$ and $Y$ be independent. Then, by Lemma \ref{lemindep}, $\langle l_{1}, X \rangle_{{\cal{H}}}$ and $\langle l_{2}, Y \rangle_{{\cal{H}}}$ are also independent, for all $l_1$ and $l_2 \in \cal{H}$. In particular, the joint probability density function $f_{{\inpr{l_1}{X}},{\inpr{l_2}{Y}}}$ equals the product of marginal probability density functions $f_{\inpr{l_1}{X}}$ and $f_{\inpr{l_2}{Y}}$. Hence, $T(R_1, R_2) = 0$.

Conversely, suppose that $T(R_1, R_2) = 0$. Then, for all $l_1$ and $l_2 \in \cal{H}$ with $\|l_1\|_{\cal{H}} \leq R_1$ and $\|l_2\|_{\cal{H}} \leq R_2$, the joint probability density function $f_{{\inpr{l_1}{X}},{\inpr{l_2}{Y}}}$ equals the product of marginal probability density functions $f_{\inpr{l_1}{X}}$ and $f_{\inpr{l_2}{Y}}$. Therefore, $\langle l_{1}, X \rangle_{{\cal{H}}}$ and $\langle l_{2}, Y \rangle_{{\cal{H}}}$ are independent.

Now, consider $\bar l_1, \bar l_2 \in \cal{H}$ with $\|\bar l_1\|_{\cal{H}} > R_1$ and $\|\bar l_2\|_{\cal{H}} > R_2$. Since $\left\|\frac{R_1}{\|\bar l_1\|_{\cal{H}}} \bar l_1\right\|_{\cal{H}} = R_1$ and $\left\|\frac{R_2}{\|\bar l_2\|_{\cal{H}}} \bar l_2\right\|_{\cal{H}} = R_2$, independence of $\langle \frac{R_1}{\|\bar l_1\|_{\cal{H}}} \bar l_1, X \rangle_{{\cal{H}}}$ and $\langle \frac{R_2}{\|\bar l_2\|_{\cal{H}}} \bar l_2, Y \rangle_{{\cal{H}}}$ follows from the hypothesis. Finally, the independence of $\langle \bar l_1, X \rangle_{{\cal{H}}}$ and $\langle \bar l_2, Y \rangle_{{\cal{H}}}$ follows from Lemma \ref{independence-under-scalar-mult}. The argument is completed by applying Lemma \ref{lemindep}.\hfill \qed

\noindent {\bf Proof of Proposition \ref{equivalence}:} Suppose that $T(R_1, R_2) = 0$. As argued in the proof of Proposition \ref{prilimarycriterionstatement}, we have for all $l_1$ and $l_2 \in \cal{H}$ with $\|l_1\|_{\cal{H}} \leq R_1$ and $\|l_2\|_{\cal{H}} \leq R_2$, the real valued random variables $\langle l_{1}, X \rangle_{{\cal{H}}}$ and $\langle l_{2}, Y \rangle_{{\cal{H}}}$ are independent.

Let $\bar l_1$ and $\bar l_2 \in \cal{H}$ with $\|\bar l_1\|_{\cal{H}} = 1$ and $\|\bar l_2\|_{\cal{H}} = 1$. Then, $\|R_1 \bar l_1\|_{\cal{H}} = R_1$ and $\|R_2\bar l_2\|_{\cal{H}} = R_2$ and consequently, the independence of $\langle R_1 \bar l_1, X \rangle_{{\cal{H}}}$ and $\langle R_2\bar l_2, Y \rangle_{{\cal{H}}}$ follows from the assertion in Proposition \ref{prilimarycriterionstatement}. Then, the independence of $\langle \bar l_1, X \rangle_{{\cal{H}}}$ and $\langle \bar l_2, Y \rangle_{{\cal{H}}}$ follows from Lemma \ref{independence-under-scalar-mult}. In particular, the joint probability density function $f_{{\inpr{\bar l_1}{X}},{\inpr{\bar l_2}{Y}}}$ equals the product of marginal probability density functions $f_{\inpr{\bar l_1}{X}}$ and $f_{\inpr{\bar l_2}{Y}}$. Since $\bar l_1$ and $\bar l_2$ are arbitrary, we have $T = 0$.

Conversely, suppose that $T = 0$, which implies that the independence of $\langle l_{1}, X \rangle_{{\cal{H}}}$ and $\langle l_{2}, Y \rangle_{{\cal{H}}}$ with $\|l_1\|_{\cal{H}} = 1$ and $\|l_2\|_{\cal{H}} = 1$. Hence, by Lemma \ref{independence-under-scalar-mult}, $\langle G_{1} l_{1}, X \rangle_{{\cal{H}}}$ and $\langle G_{2}l_{2}, Y \rangle_{{\cal{H}}}$ are independent random variables, where $0 < G_{1}\leq R_{1}$ and $0 < G_{2}\leq R_{2}$. Then, for all $l_1$ and $l_2 \in \cal{H}$ with $\|l_1\|_{\cal{H}} \leq R_1$ and $\|l_2\|_{\cal{H}} \leq R_2$, the joint probability density function $f_{{\inpr{l_1}{X}},{\inpr{l_2}{Y}}}$ equals the product of marginal probability density functions $f_{\inpr{l_1}{X}}$ and $f_{\inpr{l_2}{Y}}$. Hence, $T(R_1, R_2) = 0$, which completes the proof.  \hfill \qed

\noindent {\bf Proof of Theorem \ref{characterization}:} Note that Proposition \ref{prilimarycriterionstatement} implies that $T(R_1, R_1) = 0$ if and only if $X$ and $Y$ are independent random elements, and Proposition \ref{equivalence} implies that $T(R_1, R_2) = 0\Leftrightarrow T = 0$. Hence, both facts together imply that $T = 0$ if and only if $X$ and $Y$ are independent random elements, which completes the proof.  
\hfill \qed

\noindent {\bf Proof of Proposition \ref{lk}:} For $j = 1$ and $2$, we write $l_{j}^{(M)} = \sum\limits_{i = 1}^{M} l_{j, i} e_i$. Since, $l_j = \sum\limits_{i = 1}^{\infty} l_{j, i} e_i$ and $\|l_j\|_{\cal{H}}^2 = \sum\limits_{i = 1}^{\infty} l_{j, i}^{2} = 1$, we have
\begin{subequations}
\begin{equation}\label{apprx-by-partial-sum}
\lim_{M \to \infty}\left\|l_{j} - l_{j}^{(M)}\right\|^{2}_{\cal{H}} = \lim_{M \to \infty} \sum\limits_{i = M + 1}^{\infty} l_{j, i}^2 = 0,
\end{equation}
\text{and}
\begin{equation}\label{partial-sum-norm-bound}
\|l_{j}^{(M)}\|_{\cal{H}} \uparrow \left\|l_{j}\right\|_{\cal{H}} = 1, \text{ as } M \to \infty.
\end{equation}
\end{subequations}
Now
\begin{equation}\label{approx-by-polar-transform}
\begin{split}
\left\|l_{j} - \hat{l}_{j}^{K}\right\|_{\cal{H}} &= \left\|\left(l_{j} - l_{j}^{(M)}\right) + \left(1 - \frac{1}{\|l_{j}^{(M)}\|_{\cal{H}}}\right)l_{j}^{(M)} + \left(\frac{1}{\|l_{j}^{(M)}\|_{\cal{H}}}l_{j}^{(M)} - \hat{l}_{j}^{K}\right)\right\|_{\cal{H}}\\
&\leq \left\|l_{j} - l_{j}^{(M)}\right\|_{\cal{H}} + \left|1 - \frac{1}{\|l_{j}^{(M)}\|_{\cal{H}}}\right| \left\|l_{j}^{(M)}\right\|_{\cal{H}} + \left\|\frac{1}{\|l_{j}^{(M)}\|_{\cal{H}}}l_{j}^{(M)} - \hat{l}_{j}^{K}\right\|_{\cal{H}}\\
&\leq \left\|l_{j} - l_{j}^{(M)}\right\|_{\cal{H}} + \left|1 - \frac{1}{\|l_{j}^{(M)}\|_{\cal{H}}}\right| + \left\|\frac{1}{\|l_{j}^{(M)}\|_{\cal{H}}}l_{j}^{(M)} - \hat{l}_{j}^{K}\right\|_{\cal{H}}
\end{split}
\end{equation}
Note that
\[\left\|\frac{1}{\|l_{j}^{(M)}\|_{\cal{H}}}l_{j}^{(M)} - \hat{l}_{j}^{K}\right\|_{\cal{H}} = \left\|\left(\frac{l_{j, 1}}{\|l_{j}^{(M)}\|_{\cal{H}}}, \frac{l_{j, 2}}{\|l_{j}^{(M)}\|_{\cal{H}}}, \cdots, \frac{l_{j, M}}{\|l_{j}^{(M)}\|_{\cal{H}}}\right) - (l_{j, 1}^{(K)}, l_{j, 2}^{(K)}, \cdots, l_{j, M}^{(K)})\right\|_{\R^M},\]
with the vector $\left(\frac{l_{j, 1}}{\|l_{j}^{(M)}\|_{\cal{H}}}, \frac{l_{j, 2}}{\|l_{j}^{(M)}\|_{\cal{H}}}, \cdots, \frac{l_{j, M}}{\|l_{j}^{(M)}\|_{\cal{H}}}\right)$ in the unit sphere of $\R^M$. Therefore, by choosing large $K$ and appropriate $\hat{l}_{j}^{K}$, the term $\left\|\frac{1}{\|l_{j}^{(M)}\|_{\cal{H}}}l_{j}^{(M)} - \hat{l}_{j}^{K}\right\|_{\cal{H}}$ can be made arbitrarily small. The proof then follows by using \eqref{apprx-by-partial-sum} and \eqref{partial-sum-norm-bound} in \eqref{approx-by-polar-transform}.
\hfill \qed

\noindent The following Lemmas are useful in proving Theorem \ref{asymptotic}.

\begin{lemma}\label{ASI}
Under (A1)--(A4) and (A6), \[(I):={\sqrt{n} h_{n}}\left\{\frac{1}{n h_{n}^{{2}}}\left\{\sum\limits_{i = 1}^{n}k\left(\frac{Z_{1, i} - s}{h_{n}}, \frac{Z_{2, i} - t}{h_{n}}\right)\right\} -f_{Z_{1}, Z_{2}}(s, t)\right\}\] converges weakly to a random variable associated with Gaussian distribution with mean $= \sqrt{c}\left\{\sum\limits_{i = 1}^{2}\left(\frac{\partial f_{Z_{1}, Z_{2}}(s, t)}{\partial Z_{i}}\right)\int m_{i} k(m_1, m_2)dm_1 dm_2\right\}$ and variance $= f_{Z_1, Z_2}(s, t)\int k^{2}(m_1, m_2) dm_1 dm_2$. Here $c = \displaystyle\lim_{n\rightarrow\infty} nh_{n}^{4}$. 
\end{lemma}

\noindent {\bf Proof of Lemma \ref{ASI}:}
Observe that
\begin{eqnarray*}
&& \sqrt{n} h_{n}\left\{\frac{1}{n h_{n}^{2}}\sum\limits_{i = 1}^{n}k\left(\frac{Z_{1, i} - s}{h_{n}}, \frac{Z_{2, i} - t}{h_{n}}\right) -f_{Z_{1}, Z_{2}}(s, t)\right\}\\
& = & \sqrt{n} h_{n}\left\{\frac{1}{n h_{n}^{2}}\sum\limits_{i = 1}^{n}k\left(\frac{Z_{1, i} - s}{h_{n}}, \frac{Z_{2, i} - t}{h_{n}}\right) - E\left(\frac{1}{n h_{n}^{2}}\sum\limits_{i = 1}^{n}k\left(\frac{Z_{1, i} - s}{h_{n}}, \frac{Z_{2, i} - t}{h_{n}}\right)\right) \right\}\\
& + & \sqrt{n} h_{n}\left\{E\left(\frac{1}{n h_{n}^{2}}\sum\limits_{i = 1}^{n}k\left(\frac{Z_{1, i} - s}{h_{n}}, \frac{Z_{2, i} - t}{h_{n}}\right)\right) - f_{Z_1, Z_2} (s, t)\right\}\\
& = & (I)_{A} + (I)_{B},
\end{eqnarray*} where 
\begin{equation}\label{IA}
(I)_{A} : = \sqrt{n} h_{n}\left\{\frac{1}{n h_{n}^{2}}\sum\limits_{i = 1}^{n}k\left(\frac{Z_{1, i} - s}{h_{n}}, \frac{Z_{2, i} - t}{h_{n}}\right) - E\left(\frac{1}{n h_{n}^{2}}\sum\limits_{i = 1}^{n}k\left(\frac{Z_{1, i} - s}{h_{n}}, \frac{Z_{2, i} - t}{h_{n}}\right)\right) \right\},   
\end{equation} and 
\begin{equation}\label{IB}
(I)_{B} : =  \sqrt{n} h_{n}\left\{E\left(\frac{1}{n h_{n}^{2}}\sum\limits_{i = 1}^{n}k\left(\frac{Z_{1, i} - s}{h_{n}}, \frac{Z_{2, i} - t}{h_{n}}\right)\right) - f_{Z_1, Z_2} (s, t)\right\}.  
\end{equation}
\hfill \qed

Let us now work on \eqref{IB}.
\begin{eqnarray*}
&&\sqrt{n} h_{n}\left\{E\left(\frac{1}{n h_{n}^{2}}\sum\limits_{i = 1}^{n}k\left(\frac{Z_{1, i} - s}{h_{n}}, \frac{Z_{2, i} - t}{h_{n}}\right)\right) - f_{Z_1, Z_2} (s, t)\right\}\\
& = & \sqrt{n} h_{n}\left\{E\left(\frac{1}{h_{n}^{2}}k\left(\frac{Z_{1} - s}{h_{n}}, \frac{Z_{2} - t}{h_{n}}\right)\right) - f_{Z_1, Z_2} (s, t)\right\}\\
& = & {\sqrt{n} h_{n}\frac{1}{h_{n}^{2}}\left\{\int k\left(\frac{z_1 - s}{h_n}, \frac{z_2 - s}{h_n}\right) f_{Z_1, Z_2}(z_1, z_2) dz_1 dz_2 - h_{n}^{2}f_{Z_1, Z_2}(s, t)\right\}}\\
& = & \sqrt{n} h_{n}\left\{\int k\left(m_1, m_2\right) f_{Z_1, Z_2}(s + m_{1}h_{n}, t + m_{2}h_{n}) dm_1 dm_2 - f_{Z_1, Z_2}(s, t)\right\}\\
& = & \sqrt{n} h_{n}\left\{\int k\left(m_1, m_2\right)\left[f_{Z_1, Z_2}(s, t) + h_{n}\sum\limits_{i = 1}^{2}\left(\frac{\partial f_{Z_{1}, Z_{2}}(\xi_{n})}{\partial Z_{i}}\right) m_{i}\right]dm_{1}dm_{2} - f_{Z_1, Z_2}(s, t)\right\}\\
& = & \sqrt{n}h_{n}^{2}\left\{\sum\limits_{i = 1}^{2}\left(\frac{\partial f_{Z_{1}, Z_{2}}(\xi_{n})}{\partial Z_{i}}\right)\int m_{i} k(m_1, m_2)dm_1 dm_2 \right\}. 
\end{eqnarray*} The last fact follows from (A6), i.e. $\int k(m_1, m_2) dm_1 dm_2 = 1$. Now, using (A1), (A2) and (A6) (i.e. $h_{n}\rightarrow 0$ as $n\rightarrow\infty$, $\sqrt{n}h_{n}^{2}\rightarrow \sqrt{c}$ as $n\rightarrow\infty$, the partial derivatives of $f_{Z_1, Z_2}(., .)$ are uniformly bounded and $\int m_{i} k(m_1, m_2) dm_1 dm_2 < \infty$), we have
\begin{eqnarray}\label{LIB}
(I)_{B}\rightarrow \sqrt{c}\left\{\sum\limits_{i = 1}^{2}\left(\frac{\partial f_{Z_{1}, Z_{2}}(s, t)}{\partial Z_{i}}\right)\int m_{i} k(m_1, m_2)dm_1 dm_2 \right\}~\mbox{as $n\rightarrow\infty$}.
\end{eqnarray}

Now, we work on $(I)_{A}$ (see \eqref{IA}) and denote 
$$W_{n, i} = \sqrt{n} h_{n}\left\{\frac{1}{n h_{n}^{2}}k\left(\frac{Z_{1, i} - s}{h_{n}}, \frac{Z_{2, i} - t}{h_{n}}\right) - E\left(\frac{1}{n h_{n}^{2}}k\left(\frac{Z_{1} - s}{h_{n}}, \frac{Z_{2} - t}{h_{n}}\right)\right) \right\}, i = 1, \ldots, n$$ 
Note that $E(W_{n, i}) = 0$ for all $i = 1, \ldots, n$. Observe that 
\begin{eqnarray*}
&& W_{n, i}^{2} = \frac{1}{n h_{n}^{2}}\left\{k^{2}\left(\frac{Z_{1, i} - s}{h_{n}}, \frac{Z_{2, i} - t}{h_n}\right) + \left(E\left(k\left(\frac{Z_1 - s}{h_{n}}, \frac{Z_2 - t}{h_{n}}\right)\right)\right)^{2}\right.\\
&&\left.- 2 k\left(\frac{Z_1 - s}{h_n}, \frac{Z_2 - t}{h_n}\right) E\left(k\left(\frac{Z_1 - s}{h_n}, \frac{Z_2 - t}{h_n}\right)\right)\right\}\\
&&\Rightarrow E(W_{n, i}^{2}) = \frac{1}{n h_{n}^{2}} \left\{E\left(k^{2}\left(\frac{Z_1 - s}{h_{n}},\frac{Z_2 - t}{h_n}\right)\right) - \left(E\left(k\left(\frac{Z_1 - s}{h_{n}}, \frac{Z_2 - t}{h_{n}}\right)\right)\right)^{2}\right\}\\
& = & \frac{1}{nh_{n}^{2}}\left\{\int\left(k^{2}\left(\frac{z_1 - s}{h_{n}},\frac{z_2 - t}{h_n}\right)\right) f_{Z_1, Z_2}(z_1, z_2)dz_1 dz_2\right.\\
&&\left.- \int \left\{k\left(\frac{z_1 - s}{h_{n}},\frac{z_2 - t}{h_n}\right) f_{Z_1, Z_2}(z_1, z_2)dz_1 dz_2\right\}^{2}\right\}\\
&=&\frac{1}{n}\left\{\int k^{2}(m_1, m_2) f_{Z_1, Z_2}(s + m_1 h_{n}, t + m_2 h_{n})\right\}\\
&&- \frac{h_{n}^{2}}{n}\int \left\{k(m_1, m_2) f_{Z_1, Z_2}(s + m_1 h_n, t + m_2 h_n)\right\}^{2}dm_1 dm_2\\
&&\Rightarrow\sum\limits_{i = 1}^{n} E(W_{n, i}^{2}) = \left\{\int k^{2}(m_1, m_2) f_{Z_1, Z_2}(s + m_1 h_{n}, t + m_2 h_{n})\right\}\\
&&- h_{n}^{2}\int \left\{k(m_1, m_2) f_{Z_1, Z_2}(s + m_1 h_n, t + m_2 h_n)\right\}^{2}dm_1 dm_2.
\end{eqnarray*}
Hence, using (A1), (A2) and (A6) (i.e., $h_{n}\rightarrow 0$ as $n\rightarrow\infty$, partial derivatives of the joint density function of $(Z_1, Z_2)$ exist, and $\int k^{2}(m_1, m_2)dm_1 dm_2 < \infty$, $\int m_{1}k(m_1, m_2) dm_1 dm_2 < \infty$ and $\int m_2 k(m_1, m_2) dm_1 dm_2 < \infty$), we have 
\begin{equation}\label{LW}
\sum\limits_{i = 1}^{n} E(W_{n, i}^{2})\rightarrow f_{Z_1, Z_2}(s, t)\int k^{2}(m_1, m_2) dm_1 dm_2 : = \sigma_{1}^{2} ~\mbox{as $n\rightarrow\infty$}.    
\end{equation}
Further, it follows from the expression of $W_{n, i}$ and $(I)_{A}$ (see \eqref{IA}) that $(I)_{A} = \sum\limits_{i = 1}^{n} W_{n, i}$. Therefore, in order to prove the asymptotic normality of $(I)_{A}$, the validation of Lyapunov condition (see \cite{Serfling1980}) is established here. For some $\delta > 0$, let us consider 
\begin{eqnarray*}
&& 0\leq\frac{1}{\sigma_{1}^{2 + \delta}}\sum\limits_{i = 1}^{n}E|W_{n, i}|^{2 + \delta}\\
&&\leq\frac{1}{\sigma_{1}^{2 + \delta}}\sum\limits_{i = 1}^{n} E\left|\frac{1}{\sqrt{n}h_{n}} k\left(\frac{Z_1 - s}{h_n}, \frac{Z_2 - t}{h_n}\right)\right|^{2 + \delta}\\
&&\leq\frac{1}{\sigma_{1}^{2 + \delta}}\times\frac{1}{n^{\frac{\delta}{2}}h_{n}^{2 + \delta}}E\left|k\left(\frac{Z_1 - s}{h_n}, \frac{Z_2 - t}{h_n}\right)\right|^{2 + \delta} \leq\frac{1}{\sigma_{1}^{2 + \delta}}\times\frac{1}{n^{\frac{\delta}{2}}h_{n}^{2 + \delta}}\left|E\left(k\left(\frac{Z_1 - s}{h_n}, \frac{Z_2 - t}{h_n}\right)\right)\right|^{2 + \delta}\\
& = & \frac{1}{\sigma_{1}^{2 + \delta}}\times\frac{1}{n^{\frac{\delta}{2}}h_{n}^{2 + \delta}}\left|\int k\left(\frac{z_1 - s}{h_n}, \frac{z_2 - t}{h_n}\right) f_{Z_1, Z_2}(z_1, z_2)(z_1, z_2) dz_1 dz_2\right|^{2 + \delta}\\
& = & \frac{1}{\sigma_{1}^{2 + \delta}}\times\frac{h_{n}^{4 + 2\delta}}{n^{\frac{\delta}{2}}h_{n}^{2 + \delta}}\left|\int k\left(m_1, m_2\right) f_{Z_1, Z_2}(s + m_1 h_{n}, t + m_2 h_{n}) dm_1 dm_2\right|^{2 + \delta}\\
& = & \frac{1}{\sigma_{1}^{2 + \delta}}\times\frac{h_{n}^{2 + \delta}}{n^{\frac{\delta}{2}}}\left|\int k\left(m_1, m_2\right) f_{Z_1, Z_2}(s + m_1 h_{n}, t + m_2 h_{n}) dm_1 dm_2\right|^{2 + \delta}\rightarrow 0~\mbox{as $n\rightarrow\infty$.}
\end{eqnarray*} 
The last fact follows from the fact that $h_{n}\rightarrow 0$ as $n\rightarrow\infty$, and in view of the existence of the partial derivatives of the joint density function of $(Z_1, Z_2)$, and $\int k(m_1, m_2) dm_1 dm_2 < \infty$, $\int m_{1} k(m_1, m_2) dm_1 dm_2 <\infty$ and $\int m_2 k(m_1, m_2) dm_1 dm$. Hence, $(I)_{A}$ converges weakly to a random variable associated with Gaussian distribution with mean $= 0$ and variance $= f_{Z_1, Z_2}(s, t)\int k^{2}(m_1, m_2) dm_1 dm_2$. Further, recall that $(I) = (I)_{A} + (I)_{B}$, and it is established that $(I)_{B}\rightarrow \sqrt{c}\left\{\sum\limits_{i = 1}^{2}\left(\frac{\partial f_{Z_{1}, Z_{2}}(s, t)}{\partial Z_{i}}\right)\int m_{i} k(m_1, m_2)dm_1 dm_2 \right\}$ as $n\rightarrow\infty$ (see \eqref{LIB}), and consequently, $(I)$ converges weakly to a random variable associated with Gaussian distribution with mean $= \sqrt{c}\left\{\sum\limits_{i = 1}^{2}\left(\frac{\partial f_{Z_{1}, Z_{2}}(s, t)}{\partial Z_{i}}\right)\int m_{i} k(m_1, m_2)dm_1 dm_2 \right\}$ and variance $= f_{Z_1, Z_2}(s, t)\int k^{2}(m_1, m_2) dm_1 dm_2$.\hfill \qed

\begin{lemma}\label{ASII}
Let us denote\[(II):= \sqrt{n} h_{n}\left\{\frac{1}{n h_{n}}\sum\limits_{i = 1}^{n}k_{1}\left(\frac{Z_{1, i} - s}{h_{n}}\right)\left(\frac{1}{n h_{n}}\sum\limits_{j = 1}^{n}k_{2}\left(\frac{Z_{2, j} - t}{h_{n}}\right) - f_{Z_{2}}(t)\right)\right\}.\] Then, under (A1)--(A5),  
\begin{eqnarray}\label{6.19}
(II)\stackrel{p}\rightarrow\sqrt{c}f_{Z_{1}}(s)f^{'}_{Z_2}(t)\int mk_{2}(m) dm~\mbox{as $n\rightarrow\infty$},    
\end{eqnarray} where $c = \displaystyle\lim_{n\rightarrow\infty} nh_{n}^{4}$.
\end{lemma}

\noindent {\bf Proof of Lemma \ref{ASII}:} Note that $(II) = (II)_{A}\times (II)_{B}$, where 
\begin{equation}\label{IIA}
(II)_{A}: = \frac{1}{n h_{n}}\sum\limits_{i = 1}^{n}k_{1}\left(\frac{Z_{1, i} - s}{h_{n}}\right),   
\end{equation} and 
\begin{equation}\label{IIB}
(II)_{B}: = \sqrt{n} h_{n}\left\{\left(\frac{1}{n h_{n}}\sum\limits_{j = 1}^{n}k_{2}\left(\frac{Z_{2, j} - t}{h_{n}}\right) - f_{Z_{2}}(t)\right)\right\}.   \end{equation}
It now follows from \eqref{IIA} that 
\begin{eqnarray*}
E((II)_{A}) &=& E\left(\frac{1}{n h_{n}}\sum\limits_{i = 1}^{n}k_{1}\left(\frac{Z_{1, i} - s}{h_{n}}\right)\right)\\
& = & \frac{1}{n h_{n}}\sum\limits_{i = 1}^{n} E\left(k_{1}\left(\frac{Z_{1, i} - s}{h_{n}}\right)\right) = \frac{1}{h_{n}} E\left(k_{1}\left(\frac{Z_{1} - s}{h_{n}}\right)\right)\\
& = & \frac{1}{h_{n}}\int \left(k_{1}\left(\frac{Z_{1} - s}{h_{n}}\right)\right) f_{Z_{1}}(z_{1}) dz_{1}\\
& = & \int k(m) f_{Z_{1}} (m h_{n} + s) dm
 =  \int k(m)\left[f_{Z_{1}}(s) + mh_{n}f_{Z_{1}^{'}}(\xi_{n})\right]dm\\
 &&~\mbox{($f_{Z_{1}}^{'}$ denotes the derivative of $f_{Z_{1}}$, and $\xi_{n}\in (s, s + mh_{n})$)}\\
 &=&f_{Z_{1}}(s)\int k(m)dm + h_{n} f_{Z_{1}}^{'}(\xi_{n})\int m k(m) dm. 
\end{eqnarray*}
Now, using conditions (A1), (A4) and (A5) (i.e., $\int m k(m) dm = 1$, $h_{n}\rightarrow 0$ as $n\rightarrow\infty$, $f_{Z_{1}}^{'}(.)$ is uniformly bounded, and $\int m k(m) dm < \infty$), we have 
\begin{equation}\label{EIIA}
E((II)_{A})\rightarrow f_{Z_{1}}(s)~\mbox{as $n\rightarrow\infty$.}
\end{equation}

Next, note that 
\begin{eqnarray*}
&& 0\leq {var}\left(\frac{1}{n h_{n}}\sum\limits_{i = 1}^{n}k_{1}\left(\frac{Z_{1, i} - s}{h_{n}}\right)\right)\\
&=&\frac{1}{n^2 h_{n}^{2}}{var}\left(\sum\limits_{i = 1}^{n}k_{1}\left(\frac{Z_{1, i} - s}{h_{n}}\right)\right) + \frac{1}{n^2 h_{n}^{2}}\sum\limits_{i\neq j = 1}^{n} cov\left(k_{1}\left(\frac{Z_{1, i} - s}{h_{n}}\right), k_{1}\left(\frac{Z_{1, j} - s}{h_{n}}\right)\right) \\
& = & \frac{1}{n^2 h_{n}^{2}}var\left(\sum\limits_{i = 1}^{n}k_{1}\left(\frac{Z_{1, i} - s}{h_{n}}\right)\right)~\mbox{(since $Z_{1, i}\perp Z_{1, j}$ for $i\neq j$)}\\
& = & \frac{1}{n h_{n}^{2}} E \left(k_{1}^{2}\left(\frac{Z_1 - s}{h_{n}}\right)\right) - \frac{1}{n h_{n}^{2}} \left\{E \left(k_{1}^{2}\left(\frac{Z_1 - s}{h_{n}}\right)\right)\right\}^{2}\\
&\leq & \frac{1}{n h_{n}^{2}} E \left(k_{1}^{2}\left(\frac{Z_1 - s}{h_{n}}\right)\right)\\
& = & \frac{1}{n h_{n}^{2}}\int k_{1}^{2}\left(\frac{Z_{1} - s}{h_{n}}\right) f_{Z_{1}}(z_1) dz_1 = \frac{1}{n h_{n}}\int k_{1}^{2}(m) f_{Z_{1}} (s + m h_{n}) dm\\
& = &\frac{1}{n h_{n}}\int k_{1}^{2}(m)\left[f_{Z_{1}}(s) + m h_{n} f_{Z_{1}}^{'}(\xi_{n})\right]dm~\mbox{(here $\xi_{n}\in (s, s + m h_{n})$)}\\
& = & \frac{f_{Z_{1}}(s)}{n h_{n}}\int k_{1}^{2}(m) dm + \frac{f_{Z_{1}}^{'}(\xi_{n})}{n}\int m k_{1}^{2}(m) dm
\end{eqnarray*}. 

Now using (A1), (A3), (A4) and (A5) (i.e., $n h_{n}\rightarrow\infty$ as $n\rightarrow\infty$, $f_{Z_{1}}(.)$ and $f_{Z_{1}}^{'}(.)$ are uniformly bounded, $\int k_{1}^{2}(m) dm< \infty$ and $\int m k_{1}^{2}(m) dm < \infty$), we have 
\begin{equation}\label{VIIA}
var\left(\frac{1}{n h_{n}}\sum\limits_{i = 1}^{n}k_{1}\left(\frac{Z_{1, i} - s}{h_{n}}\right)\right)\rightarrow 0~\mbox{as $n\rightarrow\infty$}.
\end{equation}

Hence, using \eqref{EIIA} and \eqref{VIIA}, we have 
\begin{equation}\label{PIIA}
(II)_{A}\stackrel{p}\rightarrow f_{Z_{1}}(s)~\mbox{as $n\rightarrow\infty$}.    
\end{equation}

Now, let us work on $(II)_{B}$. \eqref{IIA} indicates that
\begin{eqnarray*}
(II)_{B}: &=& \sqrt{n} h_{n}\left\{\left(\frac{1}{n h_{n}}\sum\limits_{j = 1}^{n}k_{2}\left(\frac{Z_{2, j} - t}{h_{n}}\right) - f_{Z_{2}}(t)\right)\right\}\\
&=& \sqrt{n} h_{n}\left\{\left(\frac{1}{n h_{n}}\sum\limits_{j = 1}^{n}k_{2}\left(\frac{Z_{2, j} - t}{h_{n}}\right) - E\left(\frac{1}{n h_{n}}\sum\limits_{j = 1}^{n}k_{2}\left(\frac{Z_{2, j} - t}{h_{n}}\right)\right)\right)\right\}\\
& + & \sqrt{n} h_{n}\left\{E\left(\frac{1}{n h_{n}}\sum\limits_{j = 1}^{n}k_{2}\left(\frac{Z_{2, j} - t}{h_{n}}\right)\right) - f_{Z_{2}}(t)\right\}\\
& =& (II)_{B, 1} + (II)_{B, 2}, 
\end{eqnarray*}
where 
\begin{equation}\label{IIB1}
(II)_{B, 1} := \sqrt{n} h_{n}\left\{\left(\frac{1}{n h_{n}}\sum\limits_{j = 1}^{n}k_{2}\left(\frac{Z_{2, j} - t}{h_{n}}\right) - E\left(\frac{1}{n h_{n}}\sum\limits_{j = 1}^{n}k_{2}\left(\frac{Z_{2, j} - t}{h_{n}}\right)\right)\right)\right\}    
\end{equation} and 
\begin{equation}\label{IIB2}
(II)_{B, 2} :=\sqrt{n} h_{n}\left\{E\left(\frac{1}{n h_{n}}\sum\limits_{j = 1}^{n}k_{2}\left(\frac{Z_{2, j} - t}{h_{n}}\right)\right) - f_{Z_{2}}(t)\right\}.
\end{equation}

It now follows from \eqref{IIB2} that 
\begin{eqnarray*}
(II)_{B, 2} :&=&\sqrt{n} h_{n}\left\{E\left(\frac{1}{n h_{n}}\sum\limits_{j = 1}^{n}k_{2}\left(\frac{Z_{2, j} - t}{h_{n}}\right)\right) - f_{Z_{2}}(t)\right\}\\
&=& \sqrt{n} h_{n}\left\{E\left(\frac{1}{ h_{n}}k_{2}\left(\frac{Z_{2} - t}{h_{n}}\right)\right) - f_{Z_{2}}(t)\right\}\\
& = & \sqrt{n} h_{n}\left\{\frac{1}{h_{n}}\int k_{2}\left(\frac{z_2 - t}{h_{n}}\right) f_{Z_{2}}(z_2) dz_2 - f_{Z_{2}}(t)\right\}\\
& = & \sqrt{n} h_{n}\left\{\int k_{2}(m) f_{Z_{2}}(t + mh_{n}) dm - f_{Z_{2}}(t)\right\}\\
& = & \sqrt{n} h_{n}\left\{\int k_{2}(m) \left[f_{Z_{2}}(t) + m h_{n} f_{Z_{2}}^{'}(t) + \frac{m^{2}h_{n}^{2}}{2}f_{Z_{2}}^{''}(\xi_{n})\right]dm - f_{Z_{2}}(t)\right\}\\
&&~\mbox{($f_{Z_{2}}^{''}$ denotes the second derivative of $f_{Z_{2}}$), and $\xi_{n}\in (t, t + m h_{n})$}.
\end{eqnarray*}
Hence, in view of (A1), (A4) and (A5) (i.e., $h_{n}\rightarrow 0$ as $n\rightarrow\infty$, $nh_{n}^{6} = O(1)$, $\int k_{2}(m) dm = 1$, $\int m k_{2}(m) dm < \infty$, $f_{Z_{2}}^{'}$ and $f_{Z_{2}}^{''}$ are uniformly bounded), we have
\begin{equation}\label{LIIB2}
(II)_{B, 2} - \sqrt{n}h_{n}^{2}f^{'}_{Z_{2}}(t)\int mk_{2}(m)dm \rightarrow 0~\mbox{as $n\rightarrow\infty$}.
\end{equation}

Now, to study $(II)_{B, 1}$ (see \eqref{IIB1}), we have 
\begin{eqnarray}\label{6.27}
&&(II)_{B, 1} = \sqrt{n}h_{n}\left\{\left(\frac{1}{nh_{n}}\sum\limits_{j = 1}^{n}k_{2}\left(\frac{Z_{2, j} - t}{h_{n}}\right) - E\left(\frac{1}{nh_{n}}\sum\limits_{j = 1}^{n}k_{2}\left(\frac{Z_{2, j} - t}{h_{n}}\right)\right)\right)\right\}\nonumber\\
& =& \sqrt{h_{n}}\underbrace{\sqrt{nh_{n}}\left\{\left(\frac{1}{nh_{n}}\sum\limits_{j = 1}^{n}k_{2}\left(\frac{Z_{2, j} - t}{h_{n}}\right) - E\left(\frac{1}{nh_{n}}\sum\limits_{j = 1}^{n}k_{2}\left(\frac{Z_{2, j} - t}{h_{n}}\right)\right)\right)\right\}}_{\text{A}}. 
\end{eqnarray} Observe that the asymptotic normality of (A) follows using the same arguments in proving the asymptotic normality of $(I)_{A}$ (see \eqref{IA}). Therefore, as $h_{n}\rightarrow 0$ as $n\rightarrow\infty$ (see (A1)), in view of Slutsky's theorem applying on \eqref{6.27}, we have 
\begin{eqnarray}\label{6.28}
(II)_{B, 1}\stackrel{p}\rightarrow 0~\mbox{as $n\rightarrow\infty$}.   \end{eqnarray}  

As $(II)_{B} = (II)_{B, 1} + (II)_{B, 2}$, in view of \eqref{LIIB2} and \eqref{6.28} and $nh_{n}^{4}\rightarrow c$ as $n\rightarrow\infty$, one can conclude that 
\begin{eqnarray}
(II)_{B}\stackrel{p}\rightarrow \sqrt{c}f^{'}_{Z_{2}}(t)\int mk_{2}(m)dm~\mbox{as $n\rightarrow\infty$}.  \end{eqnarray} Finally, the result follows due to the fact that $(II) = (II)_{A}\times (II)_{B}$ and the assertion in \eqref{PIIA}. 
\hfill \qed

\begin{lemma}\label{ASIII}
Let us denote $$(III) : = \sqrt{n h_{n}}\left\{f_{Z_{2}}(t)\left(\frac{1}{n h_{n}}\sum\limits_{i = 1}^{n}k_{1}\left(\frac{Z_{1, i} - s}{h_{n}}\right) - f_{Z_{1}}(s)\right)\right\}.$$ Then, under (A1)--(A5), 
\begin{eqnarray}
(III)\stackrel{p}\rightarrow\sqrt{c}f_{Z_{2}}(t)f^{'}_{Z_1}(s)\int mk_{1}(m) dm~\mbox{as $n\rightarrow\infty$},    
\end{eqnarray} where $c = \displaystyle\lim_{n\rightarrow\infty} nh_{n}^{4}$.
\end{lemma}

\noindent {\bf Proof of Lemma \ref{ASIII}:} The proof of this lemma follows from the same arguments provided in the proof of Lemma \ref{ASII}. \hfill \qed

\begin{lemma}\label{pointwiseasymptotic}
Let us denote 
$$M_{n}(s, t) = \frac{1}{n h_{n}^{2}}\left\{\sum\limits_{i = 1}^{n}k\left(\frac{Z_{1, i} - s}{h_{n}}, \frac{Z_{2, i} - t}{h_{n}}\right)\right\} - \left\{\frac{1}{n^{2}h_{n}^{2}}\sum\limits_{i, j = 1}^{n}k_{1}\left(\frac{Z_{1, i} - s}{h_{n}}\right)k_{2}\left(\frac{Z_{2, j} - t}{h_{n}}\right)\right\},$$ and $$M(s, t) = f_{Z_{1}, Z_{2}}(s, t) - f_{Z_{1}}(s)f_{Z_{2}}(t),$$ where $f_{Z_{1}, Z_{2}}$, $f_{Z_{1}}$ and $f_{Z_{2}}$ are probability density functions of $(Z_1, Z_2)$, $Z_1$ and $Z_2$, respectively. Then, for any fixed $(s, t)$, under (A1)--(A6), $\sqrt{n}h_{n}(M_{n}(s, t) - M(s, t))$ converges weakly to a random variable associated with Gaussian distribution with mean $= \sqrt{c}\left\{\sum\limits_{i = 1}^{2}\left(\frac{\partial f_{Z_{1}, Z_{2}}(s, t)}{\partial Z_{i}}\right)\int m_{i} k(m_1, m_2)dm_1 dm_2\right.$\\$\left. - f_{Z_{1}}(s)f^{'}_{Z_2}(t)\int mk_{2}(m) dm - f_{Z_{2}}(t)f^{'}_{Z_1}(s)\int mk_{1}(m) dm\right\}$ \\and variance $= f_{Z_1, Z_2}(s, t)\int k^{2}(m_1, m_2) dm_1 dm_2$. Here $c = \displaystyle\lim_{n\rightarrow\infty} nh_{n}^{4}$. 
\end{lemma}

\noindent{\bf Proof of Lemma \ref{pointwiseasymptotic}:} Observe that 
\begin{eqnarray*}
&&\sqrt{n h_{n}} (M_{n} (s, t) - M(s, t))\\
& = & \sqrt{n h_{n}}\left[\frac{1}{n h_{n}^{\frac{3}{2}}}\left\{\sum\limits_{i = 1}^{n}k\left(\frac{Z_{1, i} - s}{h_{n}}, \frac{Z_{2, i} - t}{h_{n}}\right)\right\} - \left\{\frac{1}{n^2 h_{n}^2}\sum\limits_{i, j = 1}^{n}k_{1}\left(\frac{Z_{1, i} - s}{h_{n}}\right)k_{2}\left(\frac{Z_{2, j} - t}{h_{n}}\right)\right\}\right]\\
& - & \sqrt{n h_{n}} \left(f_{Z_1, Z_2}(s, t) - f_{Z_1}(s) f_{Z_2} (t)\right)\\
&=& \sqrt{n h_{n}}\left\{\frac{1}{n h_{n}^{\frac{3}{2}}}\left\{\sum\limits_{i = 1}^{n}k\left(\frac{Z_{1, i} - s}{h_{n}}, \frac{Z_{2, i} - t}{h_{n}}\right)\right\} -f_{Z_{1}, Z_{2}}(s, t)\right\}\\
& - & \sqrt{nh_{n}}\left\{\frac{1}{n^{2} h_{n}^{2}} \sum\limits_{i, j = 1}^{n}k_{1}\left(\frac{Z_{1, i} - s}{h_{n}}\right)k_{2}\left(\frac{Z_{2, j} - t}{h_{n}}\right) - f_{Z_{1}}(s) f_{Z_{2}}(t) \right\}\\
& = & \sqrt{n h_{n}}\left\{\frac{1}{n h_{n}^{\frac{3}{2}}}\left\{\sum\limits_{i = 1}^{n}k\left(\frac{Z_{1, i} - s}{h_{n}}, \frac{Z_{2, i} - t}{h_{n}}\right)\right\} -f_{Z_{1}, Z_{2}}(s, t)\right\}\\
&-& \sqrt{n h_{n}}\left\{\frac{1}{n h_{n}}\sum\limits_{i = 1}^{n}k_{1}\left(\frac{Z_{1, i} - s}{h_{n}}\right)\left(\frac{1}{n h_{n}}\sum\limits_{j = 1}^{n}k_{2}\left(\frac{Z_{2, j} - t}{h_{n}}\right) - f_{Z_{2}}(t)\right)\right\}\\
&-& \sqrt{n h_{n}}\left\{f_{Z_{2}}(t)\left(\frac{1}{n h_{n}}\sum\limits_{i = 1}^{n}k_{2}\left(\frac{Z_{1, i} - s}{h_{n}}\right) - f_{Z_{1}}(s)\right)\right\}\\
&=& (I) - (II) - (III), 
\end{eqnarray*} where 
\begin{equation}\label{I}
(I) : = \sqrt{n h_{n}}\left\{\frac{1}{n h_{n}^{\frac{3}{2}}}\left\{\sum\limits_{i = 1}^{n}k\left(\frac{Z_{1, i} - s}{h_{n}}, \frac{Z_{2, i} - t}{h_{n}}\right)\right\} -f_{Z_{1}, Z_{2}}(s, t)\right\},     
\end{equation}
\begin{equation}\label{II}
(II) : = \sqrt{n h_{n}}\left\{\frac{1}{n h_{n}}\sum\limits_{i = 1}^{n}k_{1}\left(\frac{Z_{1, i} - s}{h_{n}}\right)\left(\frac{1}{n h_{n}}\sum\limits_{j = 1}^{n}k_{2}\left(\frac{Z_{2, j} - t}{h_{n}}\right) - f_{Z_{2}}(t)\right)\right\}, \end{equation} and 
\begin{equation}\label{III}
(III) : = \sqrt{n h_{n}}\left\{f_{Z_{2}}(t)\left(\frac{1}{n h_{n}}\sum\limits_{i = 1}^{n}k_{2}\left(\frac{Z_{1, i} - s}{h_{n}}\right) - f_{Z_{1}}(s)\right)\right\}.
\end{equation}

Now, the result directly follows from the assertions in Lemmas \ref{ASI}, \ref{ASII} and \ref{ASIII} using  Slutsky's theorem. 

\hfill \qed

\begin{lemma}\label{rate-difference-T-TGL}
Under (A1), (A5) and (A6), $\sqrt{n}h_{n}\left\{\displaystyle\lim_{G\rightarrow\infty}\displaystyle\lim_{L\rightarrow\infty} (T_{n} - T_{n}^{G, L})\right\}\stackrel{p}\rightarrow 0$ as $n\rightarrow\infty$.
\end{lemma}

\noindent{\bf Proof of Lemma \ref{rate-difference-T-TGL}:}
Recall $T_{n}$ from \eqref{statistic} : 
\begin{eqnarray*}T_{n} &=& 
\sup_{\substack{\theta_{j, i}\in (-\frac{\pi}{2}, \frac{\pi}{2})\\ j = 1, 2; i = 1, \ldots, M - 2\\ \theta_{j, M - 1}\in (-\pi, \pi)\\ j = 1, 2\\ s\in\mathbb{R}, t\in\mathbb{R}}}
\left|\frac{1}{nh_{n}^{2}}\sum\limits_{i = 1}^{n}k\left(\frac{\langle \hat{l}_{1}^{K}, X_{i}\rangle_{\cal{H}} - s}{h_{n}}, \frac{\langle \hat{l}_{2}^{K}, Y_{i}\rangle_{\cal{H}} - t}{h_{n}}\right)\right.\\ - &&\left.\frac{1}{n^{2}h_{n}^{2}}\sum\limits_{i, j = 1}^{n}k_{1}\left(\frac{\langle \hat{l}_{1}^{K}, X_{i}\rangle_{\cal{H}} - s}{h_{n}}\right)k_{2}\left(\frac{\langle \hat{l}_{2}^{K}, Y_{j}\rangle_{\cal{H}} - t}{h_{n}}\right)\right|,\end{eqnarray*} and recall from
\eqref{practicestatistic} that 
\begin{eqnarray*}
T_{n}^{G, L} &=& 
\sup_{\substack{\theta_{j, i}\in (-\frac{\pi}{2}, \frac{\pi}{2})\\ j = 1, 2; i = 1, \ldots, M - 1\\ \theta_{j, M - 2}\in (-\pi, \pi)\\ j = 1, 2\\ s\in\{s_1, \ldots, s_L\}\\ t\in\{t_1, \ldots, t_L\}}}
\left|\frac{1}{n h_{n}^{2}}\sum\limits_{i = 1}^{n}k\left(\frac{\langle \hat{l}_{1}^{K}, X_{i}\rangle_{\cal{H}} - s}{h_{n}}, \frac{\langle \hat{l}_{2}^{K}, Y_{i}\rangle_{\cal{H}} - t}{h_{n}}\right)\right.\\ - &&\left.\frac{1}{n^{2}h_{n}^{2}}\sum\limits_{i, j = 1}^{n}k_{1}\left(\frac{\langle \hat{l}_{1}^{K}, X_{i}\rangle_{\cal{H}} - s}{h_{n}}\right)k_{2}\left(\frac{\langle \hat{l}_{2}^{K}, Y_{j}\rangle_{\cal{H}} - t}{h_{n}}\right)\right|.\end{eqnarray*} 

Note that it follows from the assertion in Proposition \ref{GL} that $$\sup_{\substack{s\in\{s_1, \ldots, s_L\}\\ t\in\{t_1, \ldots, t_L\}}}\sum\limits_{i = 1}^{n}k\left(\frac{\langle \hat{l}_{1}^{K}, X_{i}\rangle_{\cal{H}} - s}{h_{n}}, \frac{\langle \hat{l}_{2}^{K}, Y_{i}\rangle_{\cal{H}} - t}{h_{n}}\right) - \sup_{\substack{s\in\mathbb{R}\\ t\in\mathbb{R}}}\sum\limits_{i = 1}^{n}k\left(\frac{\langle \hat{l}_{1}^{K}, X_{i}\rangle_{\cal{H}} - s}{h_{n}}, \frac{\langle \hat{l}_{2}^{K}, Y_{i}\rangle_{\cal{H}} - t}{h_{n}}\right)\rightarrow 0$$ almost surely as $G\rightarrow\infty$ and $L\rightarrow\infty$ (iteratively) for any $n\geq 1$, where $s_{i}\in [-G, G]$ and $t_{i}\in [-G, G]$ for $i = 1, \ldots, L$. This implies that 
\begin{equation}
E\left[\sup_{\substack{s\in\{s_1, \ldots, s_L\}\\ t\in\{t_1, \ldots, t_L\}}}\sum\limits_{i = 1}^{n}k\left(\frac{\langle \hat{l}_{1}^{K}, X_{i}\rangle_{\cal{H}} - s}{h_{n}}, \frac{\langle \hat{l}_{2}^{K}, Y_{i}\rangle_{\cal{H}} - t}{h_{n}}\right) - \sup_{\substack{s\in\mathbb{R}\\ t\in\mathbb{R}}}\sum\limits_{i = 1}^{n}k\left(\frac{\langle \hat{l}_{1}^{K}, X_{i}\rangle_{\cal{H}} - s}{h_{n}}, \frac{\langle \hat{l}_{2}^{K}, Y_{i}\rangle_{\cal{H}} - t}{h_{n}}\right)\right] \rightarrow 0    
\end{equation} for any $n\geq 1$ as $G\rightarrow\infty$ and $L\rightarrow\infty$ (iteratively), which follows from the dominated convergence theorem (see, e.g., \cite{Serfling1980}) because $$\left[\sup_{\substack{s\in\{s_1, \ldots, s_L\}\\ t\in\{t_1, \ldots, t_L\}}}\sum\limits_{i = 1}^{n}k\left(\frac{\langle \hat{l}_{1}^{K}, X_{i}\rangle_{\cal{H}} - s}{h_{n}}, \frac{\langle \hat{l}_{2}^{K}, Y_{i}\rangle_{\cal{H}} - t}{h_{n}}\right) - \sup_{\substack{s\in\mathbb{R}\\ t\in\mathbb{R}}}\sum\limits_{i = 1}^{n}k\left(\frac{\langle \hat{l}_{1}^{K}, X_{i}\rangle_{\cal{H}} - s}{h_{n}}, \frac{\langle \hat{l}_{2}^{K}, Y_{i}\rangle_{\cal{H}} - t}{h_{n}}\right)\right]$$ is bounded for all $G$ and $L$ and any $n\geq 1$.  Hence, 
\begin{equation}\label{Econ}
\begin{aligned}
&E\left[\sup_{\substack{s\in\{s_1, \ldots, s_L\}\\ t\in\{t_1, \ldots, t_L\}}}\frac{\sqrt{n}h_{n}}{nh_{n}^{2}}\sum\limits_{i = 1}^{n}k\left(\frac{\langle \hat{l}_{1}^{K}, X_{i}\rangle_{\cal{H}} - s}{h_{n}}, \frac{\langle \hat{l}_{2}^{K}, Y_{i}\rangle_{\cal{H}} - t}{h_{n}}\right)\right.\\ &\left.- \sup_{\substack{s\in\mathbb{R}\\ t\in\mathbb{R}}}\frac{\sqrt{n}h_{n}}{nh_{n}^{2}}\sum\limits_{i = 1}^{n}k\left(\frac{\langle \hat{l}_{1}^{K}, X_{i}\rangle_{\cal{H}} - s}{h_{n}}, \frac{\langle \hat{l}_{2}^{K}, Y_{i}\rangle_{\cal{H}} - t}{h_{n}}\right)\right] \rightarrow 0 
\end{aligned}
\end{equation} as $n\rightarrow\infty$ (along with $G\rightarrow\infty$ and $L\rightarrow\infty$, iteratively) in view of (A1).

\noindent Arguing exactly in a similar way, one can conclude that 
\begin{equation}\label{Vcon}
\begin{aligned}
 &E\left[\sup_{\substack{s\in\{s_1, \ldots, s_L\}\\ t\in\{t_1, \ldots, t_L\}}}\frac{\sqrt{n}h_{n}}{nh_{n}^{2}}\sum\limits_{i = 1}^{n}k\left(\frac{\langle \hat{l}_{1}^{K}, X_{i}\rangle_{\cal{H}} - s}{h_{n}}, \frac{\langle \hat{l}_{2}^{K}, Y_{i}\rangle_{\cal{H}} - t}{h_{n}}\right)\right.\\ &\left.- \sup_{\substack{s\in\mathbb{R}\\ t\in\mathbb{R}}}\frac{\sqrt{n}h_{n}}{nh_{n}^{2}}\sum\limits_{i = 1}^{n}k\left(\frac{\langle \hat{l}_{1}^{K}, X_{i}\rangle_{\cal{H}} - s}{h_{n}}, \frac{\langle \hat{l}_{2}^{K}, Y_{i}\rangle_{\cal{H}} - t}{h_{n}}\right)\right]^{2}\rightarrow 0 
 \end{aligned}
\end{equation} as $n\rightarrow\infty$ (along with $G\rightarrow\infty$ and $L\rightarrow\infty$ (iteratively)).

\noindent Hence, using \eqref{Econ} and \eqref{Vcon}, we have
\begin{equation}\label{jointcon}
\begin{aligned}
&\left\{\sup_{\substack{s\in\{s_1, \ldots, s_L\}\\ t\in\{t_1, \ldots, t_L\}}}\frac{\sqrt{n}h_{n}}{nh_{n}^{2}}\sum\limits_{i = 1}^{n}k\left(\frac{\langle \hat{l}_{1}^{K}, X_{i}\rangle_{\cal{H}} - s}{h_{n}}, \frac{\langle \hat{l}_{2}^{K}, Y_{i}\rangle_{\cal{H}} - t}{h_{n}}\right)\right.\\
&\left.- \sup_{\substack{s\in\mathbb{R}\\ t\in\mathbb{R}}}\frac{\sqrt{n}h_{n}}{nh_{n}^{2}}\sum\limits_{i = 1}^{n}k\left(\frac{\langle \hat{l}_{1}^{K}, X_{i}\rangle_{\cal{H}} - s}{h_{n}}, \frac{\langle \hat{l}_{2}^{K}, Y_{i}\rangle_{\cal{H}} - t}{h_{n}}\right)\right\}\stackrel{p}\rightarrow 0
\end{aligned}
\end{equation} as $n\rightarrow\infty$ (along with $G\rightarrow\infty$ and $L\rightarrow\infty$ (iteratively)).

\noindent Similarly, one can show that 
\begin{equation}\label{margcon}
\begin{aligned}
&\frac{\sqrt{n}h_{n}}{n^{2}h_{n}^{2}}\left\{\sup_{\substack{s\in\{s_1, \ldots, s_L\}\\ t\in\{t_1, \ldots, t_L\}}}\sum\limits_{i = 1}^{n}k_{1}\left(\frac{\langle \hat{l}_{1}^{K}, X_{i}\rangle_{\cal{H}} - s}{h_{n}}\right) k_{2}\left(\frac{\langle \hat{l}_{2}^{K}, Y_{i}\rangle_{\cal{H}} - t}{h_{n}}\right)\right.\\ &\left.- \sup_{\substack{s\in\mathbb{R}\\ t\in\mathbb{R}}}\sum\limits_{i = 1}^{n}k_{1}\left(\frac{\langle \hat{l}_{1}^{K}, X_{i}\rangle_{\cal{H}} - s}{h_{n}}\right) k_{2}\left(\frac{\langle \hat{l}_{2}^{K}, Y_{i}\rangle_{\cal{H}} - t}{h_{n}}\right)\right\}\stackrel{p}\rightarrow 0. \end{aligned}  
\end{equation}

\noindent In view of \eqref{jointcon} and \eqref{margcon}, we have $\sqrt{n}h_{n}\left\{\displaystyle\lim_{G\rightarrow\infty}\displaystyle\lim_{L\rightarrow\infty} (T_{n} - T_{n}^{G, L})\right\}\stackrel{p}\rightarrow 0$ as $n\rightarrow\infty$, which completes the proof. \hfill \qed

\begin{lemma}\label{pointwiseasymptotic-multidim}
For distinct pairs $(s_1, t_1), \cdots, (s_m, t_m)$, consider the random vector
\[\sqrt{nh_n}\begin{pmatrix}
M_n(s_1, t_1) - M(s_1, t_1)\\
M_n(s_2, t_2) - M(s_2, t_2)\\
\cdots\\
M_n(s_m, t_m) - M(s_m, t_m)
\end{pmatrix}
\]
where $M_n(\cdot, \cdot)$ and $M(\cdot, \cdot)$ are the same as in Lemma \ref{pointwiseasymptotic}. Then, under (A1)--(A6), the sequence of random vectors converge in distribution to an $m$-dimensional multivariate Normal distribution with independent components, where the $l$-th component follows a univariate Normal distribution with mean\\ $=\sqrt{c}\left\{\sum\limits_{i = 1}^{2}\left(\frac{\partial f_{Z_{1}, Z_{2}}(s_l, t_l)}{\partial Z_{i}}\right)\int m_{i} k(m_1, m_2)dm_1 dm_2 \right.$\\$\left. - f_{Z_{1}}(s_l)f^{'}_{Z_2}(t_l)\int mk_{2}(m) dm - f_{Z_{2}}(t_l)f^{'}_{Z_1}(s_l)\int mk_{1}(m) dm\right\}$\\ and variance $= f_{Z_1, Z_2}(s_l, t_l)\int k^{2}(m_1, m_2) dm_1 dm_2$. Here $c = \displaystyle\lim_{n\rightarrow\infty} nh_{n}^{4}.$
\end{lemma}

\noindent {\bf Proof of Lemma \ref{pointwiseasymptotic-multidim}:}
We proceed as in Lemma \ref{pointwiseasymptotic}. Here, we fix $\alpha_1, \cdots. \alpha_m \in \R$ and look at the convergence in distribution of $\sqrt{nh_n}\sum\limits_{l = 1}^m \alpha_l (M_n(s_l, t_l) - M(s_l, t_l))$ as $n$ goes to $\infty$. Observe that 
\begin{equation}
\sqrt{nh_n}\sum_{l = 1}^m \alpha_l (M_n(s_l, t_l) - M(s_l, t_l)) = \sum_{l = 1}^m \alpha_l . (I)_{s_l, t_l} - \sum_{l = 1}^m \alpha_l . (II)_{s_l, t_l} - \sum_{l = 1}^m \alpha_l . (III)_{s_l, t_l},
\end{equation}
where
\[(I)_{s_l, t_l} = \sqrt{n h_{n}}\left\{\frac{1}{n h_{n}^{\frac{3}{2}}}\left\{\sum\limits_{i = 1}^{n}k\left(\frac{Z_{1, i} - s_l}{h_{n}}, \frac{Z_{2, i} - t_l}{h_{n}}\right)\right\} -f_{Z_{1}, Z_{2}}(s_l, t_l)\right\},\]
\[(II)_{s_l, t_l} = \sqrt{n h_{n}}\left\{\frac{1}{n h_{n}}\sum\limits_{i = 1}^{n}k_{1}\left(\frac{Z_{1, i} - s_l}{h_{n}}\right)\left(\frac{1}{n h_{n}}\sum\limits_{j = 1}^{n}k_{2}\left(\frac{Z_{2, j} - t_l}{h_{n}}\right) - f_{Z_{2}}(t_l)\right)\right\}\]
and
\[(III)_{s_l, t_l} = \sqrt{n h_{n}}\left\{f_{Z_{2}}(t_l)\left(\frac{1}{n h_{n}}\sum\limits_{i = 1}^{n}k_{2}\left(\frac{Z_{1, i} - s_l}{h_{n}}\right) - f_{Z_{1}}(s_l)\right)\right\}\]
for $l = 1, \cdots, m$.

Now, the result directly follows from the assertions in Lemma \ref{ASI}, \ref{ASII} and \ref{ASIII} using  Slutsky's theorem in view of the fact that $m$ is a fixed integer (independent of $n$), and for all $i\in \{1, \ldots, m\}$, $\alpha_{i}$ is a fixed constant. \hfill \qed

\noindent {\bf Proof of Theorem \ref{asymptotic}:} The proof follows from the proof of Lemma \ref{rate-difference-T-TGL} and the proof of Lemma \ref{pointwiseasymptotic-multidim} along with an application of continuous mapping theorem.
\hfill \qed

\noindent {\bf Proof of Proposition \ref{consistency}:} 
Observe that the asymptotic power of the test under $H_{1}$ is 
\begin{eqnarray*}
&&P_{H_{1}}[\sqrt{n} h_{n}T_{n}^{G, L}\geq\hat{c}_{\alpha}]~\mbox{(since under $H_{0}$, $T = 0$)}\\
&&= P[\sqrt{n} h_{n}(T_{n}^{G, L} - t_{0})\geq\hat{c}_{\alpha} - \sqrt{n }h_{n}t_{0}]~\mbox{(assume that $T = t_{0} > 0$ under $H_{1}$)}\\
&&= P[Z\geq\hat{c}_{\alpha} - \sqrt{n} h_{n}t_{0}],
\end{eqnarray*} where $Z$ follows a certain non-degenerate distribution. Since $\sqrt{n} h_{n}\rightarrow\infty$ as $n\rightarrow\infty$ (see condition (A1)) and $t_{0} > 0$, we have $$\lim_{n\rightarrow\infty}P[Z\geq\hat{c}_{\alpha} - \sqrt{n} h_{n}t_{0}] = 1.$$ This completes that proof. 
\hfill \qed

The following lemma is used in the proof of Proposition \ref{GL}.

\begin{lemma}\label{from-max-on-grid-to-sup}
Let $f:\R^2 \to \R$ be any continuous function. Let $\|f\|_\infty$ denote $\displaystyle\sup_{(t, s) \in \R^2}|f(t, s)|$ and $\|f\|_{\infty, G}$ denote $\displaystyle\sup_{(t, s) \in [-G, G] \times [-G, G]}|f(t, s)|$, for any $G > 0$. Given any $G > 0$ and any positive integer $L$, consider the set 
\[S_{G, L} := \{-G + i \tfrac{2G}{L}: i = 0, 1, \cdots, L\}\]
of equally spaced points in the interval $[-G, G]$. Then, we have the following limits.
\begin{enumerate}[label=(\roman*)]
\item For any $G > 0$,
\[\lim_{L \to \infty} \max_{(t, s) \in S_{G, L} \times S_{G, L}} |f(t, s)| = \|f\|_{\infty, G}.\]
\item We have \[\lim_{G \to \infty} \lim_{L \to \infty} \max_{(t, s) \in S_{G, L} \times S_{G, L}} |f(t, s)| = \lim_{G \to \infty} \|f\|_{\infty, G} = \|f\|_\infty.\]
\end{enumerate}
\end{lemma}

\noindent {\bf Proof of Lemma \ref{from-max-on-grid-to-sup}:}
Fix $G > 0$ and note that the continuous function $f$ is bounded on $[-G, G] \times [-G, G]$ (see \cite[Theorem 4.15]{Rudin1976}). Further, since the continuous function $|f|$ on $[-G, G] \times [-G, G]$ achieves its supremum within the set (see \cite[Theorem 4.16]{Rudin1976}), there exists $(t^\prime, s^\prime) \in [-G, G] \times [-G, G]$ such that $|f(t^\prime, s^\prime)| = \|f\|_{\infty, G}$.\\
Since $f$ is continuous at $(t^\prime, s^\prime)$, corresponding to the above $\epsilon > 0$, we have $\delta > 0$ such that 
\[|f(t, s) - f(t^\prime, s^\prime)| < \tfrac{\epsilon}{2}\]
whenever $(t, s) \in [-G, G] \times [-G, G]$ with $\sqrt{(t - t^\prime)^2 + (s - s^\prime)^2} < \delta$.\\
Given a positive integer $L$, we can cover the set $[-G, G] \times [-G, G]$ by squares with corners from $S_{G, L} \times S_{G, L}$ and with lengths of a side $\tfrac{2G}{L}$. Here, the diagonals in these squares are of length $\tfrac{2\sqrt{2}G}{L}$, and as such, any point in such a square is within a distance $\tfrac{\sqrt{2}G}{L}$ from a corner. In particular, for $(t^\prime, s^\prime)$, there exists $(t^{\prime\prime}, s^{\prime\prime}) \in S_{G, L} \times S_{G, L}$ such that $\sqrt{(t^{\prime\prime} - t^\prime)^2 + (s^{\prime\prime} - s^\prime)^2} \leq \tfrac{\sqrt{2}G}{L}$. We can now choose large $L$ such that $\tfrac{\sqrt{2}G}{L} < \delta$ and hence, $|f(t^{\prime\prime}, s^{\prime\prime}) - f(t^\prime, s^\prime)| < \tfrac{\epsilon}{2}$ and in particular,
\[|f(t^\prime, s^\prime)| < |f(t^{\prime\prime}, s^{\prime\prime})| + \tfrac{\epsilon}{2}.\]
Then,
\[|f(t^{\prime\prime}, s^{\prime\prime})| \leq \|f\|_{\infty, G} < |f(t^{\prime\prime}, s^{\prime\prime})| + \epsilon,\]
with $(t^{\prime\prime}, s^{\prime\prime}) \in S_{G, L} \times S_{G, L}$ for sufficiently large $L$.\\
Since $\epsilon > 0$ is arbitrary, for any fixed $G > 0$, we have proved the first statement, i.e.,
\[\lim_{L \to \infty} \max_{(t, s) \in S_{G, L} \times S_{G, L}} |f(t, s)| = \|f\|_{\infty, G}.\]

To prove the second statement, it is enough to show that $\lim_{G \to \infty} \|f\|_{\infty, G} = \|f\|_\infty$.

First, consider the case when $\|f\|_\infty < \infty$. Then, given $\epsilon > 0$, there exists $(t^\prime, s^\prime) \in \R^2$ such that \[|f(t^\prime, s^\prime)| \leq \|f\|_\infty < |f(t^\prime, s^\prime)| + \epsilon.\]
Then, there exists $G > 0$ large such that $(t^\prime, s^\prime) \in [-G, G] \times [-G, G]$ and hence
\[|f(t^\prime, s^\prime)| \leq \|f\|_{\infty, G} \leq \|f\|_\infty < |f(t^\prime, s^\prime)| + \epsilon \leq \|f\|_{\infty, G} + \epsilon.\]
Hence, we have $\lim_{G \to \infty} \|f\|_{\infty, G} = \|f\|_\infty$ when $\|f\|_\infty < \infty$.

When $\|f\|_\infty = \infty$, for any $R >0$, there exists $(t^\prime, s^\prime) \in \R^2$ such that \[|f(t^\prime, s^\prime)| \geq R.\]
Then, there exists $G > 0$ large such that $(t^\prime, s^\prime) \in [-G, G] \times [-G, G]$ and hence
\[R \leq |f(t^\prime, s^\prime)| \leq \|f\|_{\infty, G}.\]
Hence, we have $\lim_{G \to \infty} \|f\|_{\infty, G} = \infty = \|f\|_\infty$. This completes the proof.
\hfill \qed

\noindent {\bf Proof of Proposition \ref{GL}:} First, we observe that
given finitely many real-valued continuous functions $f_1, f_2, \cdots, f_m$ on $\R^2$, the function $(t, s) \mapsto \max\{f_1(t, s), f_2(t, s), \cdots, f_m(t, s)\}$ is also continuous. To see this, note that $(t, s) \mapsto \max\{f_1(t, s), f_2(t, s)\} = \tfrac{1}{2}|f_1(t, s) + f_2(t, s)| - \tfrac{1}{2}(f_1(t, s) - f_2(t, s))$ is a continuous function on $\R^2$. Consequently, $(t, s) \mapsto \max\{f_1(t, s), f_2(t, s), f_3(t, s)\} = \max\{ \max\{f_1(t, s), f_2(t, s)\}, f_3(t, s)\}$ is also continuous. Iterating this way, we have the above observation.

Now, for every fixed $n \geq 1$, consider the real valued continuous function
\begin{equation*}
\begin{aligned}
&(t, s) \mapsto \max_{\substack{\theta_{j, i}\in (-\frac{\pi}{2}, \frac{\pi}{2})\\ j = 1, 2; i = 1, \ldots, M - 1\\ \theta_{j, M - 2}\in (-\pi, \pi)\\ j = 1, 2}}
\left|\frac{1}{n h_{n}^{2}}\sum\limits_{i = 1}^{n}k\left(\frac{\langle \hat{l}_{1}^{K}, X_{i}\rangle_{\cal{H}} - s}{h_{n}}, \frac{\langle \hat{l}_{2}^{K}, Y_{i}\rangle_{\cal{H}} - t}{h_{n}}\right)\right.\\ - &\left.\frac{1}{n^{2}h_{n}^{2}}\sum\limits_{i, j = 1}^{n}k_{1}\left(\frac{\langle \hat{l}_{1}^{K}, X_{i}\rangle_{\cal{H}} - s}{h_{n}}\right)k_{2}\left(\frac{\langle \hat{l}_{2}^{K}, Y_{j}\rangle_{\cal{H}} - t}{h_{n}}\right)\right|
\end{aligned}
\end{equation*}
on $\R^2$. Under (A6), the required almost sure convergence follows from Lemma \ref{from-max-on-grid-to-sup}. 
\hfill \qed

As a consequence of assumption (A3) that the random elements $X$ and $Y$ are norm bounded and  the joint probability density function of the random variables $\langle l_{1}, X \rangle_{\cal{H}}$ and $\langle l_{2}, Y \rangle_{\cal{H}}$ is uniformly bounded, we show that the marginal probability density functions of the random variables $\langle l_{1}, X \rangle_{\cal{H}}$ and $\langle l_{2}, Y \rangle_{\cal{H}}$ are also bounded.
\begin{proposition}\label{propbounded}
Suppose that the Hilbert valued random elements $X$ and $Y$ are essentially bounded, i.e., there exist $M_1 > 0$ and $M_2 > 0$ such that $\p(\|X\|_{\cal{H}} \leq M_1) = 1$ and $\p(\|Y\|_{\cal{H}} \leq M_2) = 1$. Then, for all $l_1$ and $l_2 \in \cal{H}$, the random variables $\inpr{l_1}{X}$ and $\inpr{l_2}{Y}$ are also essentially bounded. Furthermore, in this case, if the joint probability density function $f_{\inpr{l_1}{X}, \inpr{l_2}{Y}}$ is bounded, then so are the marginal probability density functions $f_{\inpr{l_1}{X}}$ and $f_{\inpr{l_2}{Y}}$.
\end{proposition}

\noindent {\bf Proof:}
For $x, y \in \cal{H}$, we have $|\inpr{x}{y}| \leq \|x\|_{\cal{H}} \|y\|_{\cal{H}}$. Then, for all $l_1, l_2 \in \cal{H}$, $\p(|\inpr{l_1}{X}| \leq \|l_1\|_{\cal{H}} M_1) = 1$ and $\p(|\inpr{l_2}{Y}| \leq \|l_2\|_{\cal{H}} M_2) = 1$. This proves that the random variables $\inpr{l_1}{X}$ and $\inpr{l_2}{Y}$ are essentially bounded.

Now, if $f_{\inpr{l_1}{X}, \inpr{l_2}{Y}}(x, y) \leq R$ for all $x \in [-M_1, M_1], y \in [-M_2, M_2]$, then \[f_{\inpr{l_1}{X}}(x) = \int_{-M_2}^{M_2} f_{\inpr{l_1}{X}, \inpr{l_2}{Y}}(x, t)\, dt \leq 2RM_2, \forall x\]
and 
\[f_{\inpr{l_2}{Y}}(y) = \int_{-M_1}^{M_1} f_{\inpr{l_1}{X}, \inpr{l_2}{Y}}(s, y)\, ds \leq 2RM_1, \forall y\]
This completes the proof. \hfill \qed

\section{Appendix B : Additional Real Data Analysis and Simulation Studies}\label{ANS}
\subsection{Real Data Analysis}
\begin{figure}
	\begin{center}
	\includegraphics[width= 4.5in, height = 3.25in]{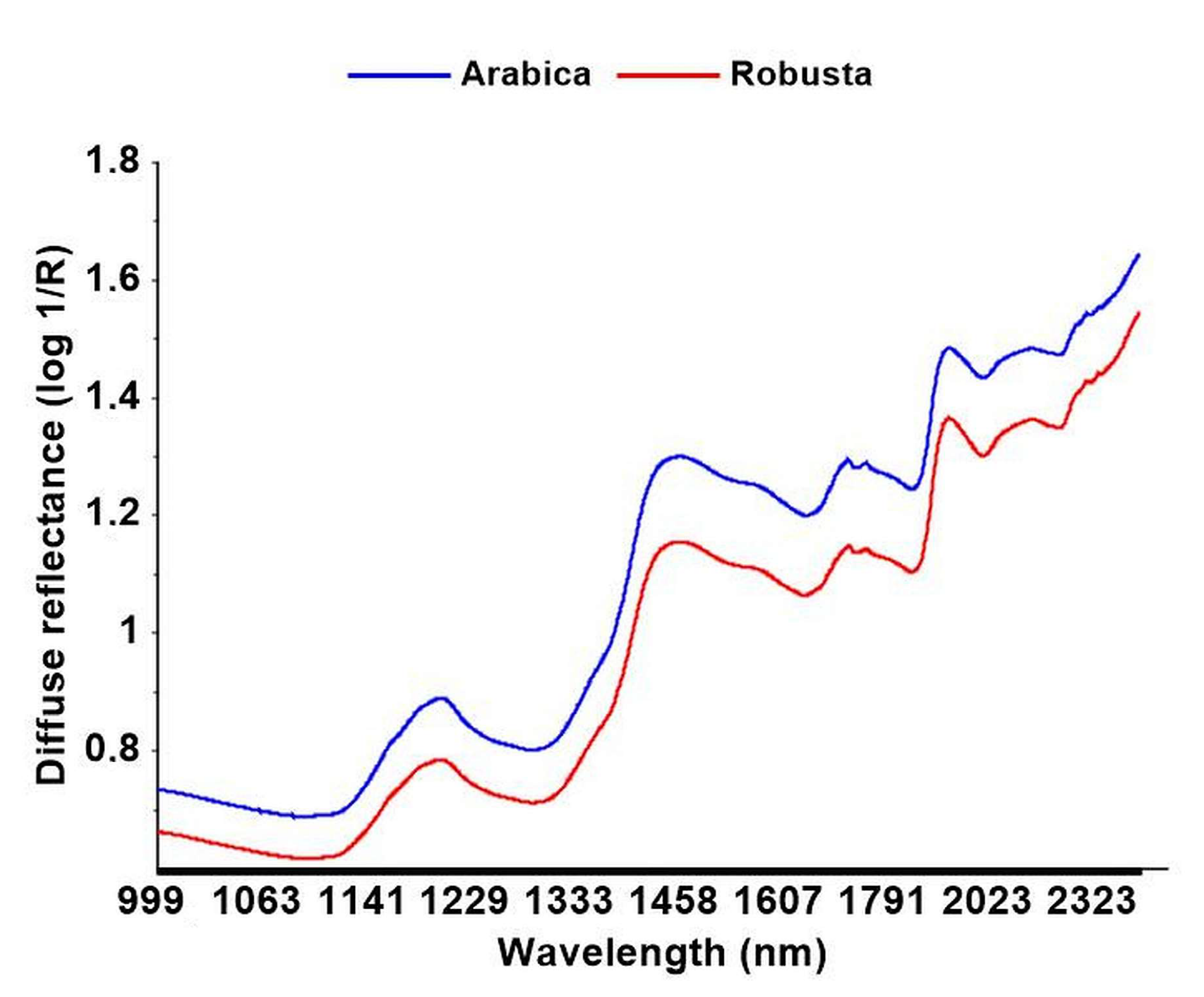}
	\end{center}
	\caption{\it  
		\label{Coffee}
		Spectroscopy readings of mean curve of Arabica coffee and the Robusta coffee. Source : \url{https://www.mdpi.com/2304-8158/9/6/788/htm}
	} 
\end{figure}

 In this section, we analyse two more real data sets. These data sets are well-known as the Coffee data and the Berkeley growth data. The coffee data is available at \url{https://www.cs.ucr.edu/~eamonn/time_series_data_2018/}, and it contains spectroscopy readings taken at 286 wavelength values for 14 samples of each of the two varieties of coffee, namely, Arabica and Robusta. The readings are illustrated in the diagrams in Figure \ref{Coffee}, and each observation is viewed as an element  in the separable Hilbert space $L_{2}([0, 1])$ after normalization by an appropriate scaling factor as we did in the analysis of Canadian weather data in Section \ref{RDA}. Following the same methodology described in Section \ref{RDA}, we compute the $p$-value, and the $p$-value is obtained as $0.069$, which indicates that the data does not favour the null hypothesis at 7\% level of significance, i.e., spectroscopy readings of the Arabica coffee and the Robusta coffee have strong enough statistical dependence. Note that this result is not unexpected as it seems from the diagrams in Figure \ref{Coffee} that the curves associated with spectroscopy readings of the Arabica coffee and the Robusta coffee are not independent. Now, as this data is not favouring the null hypothesis, here also, we estimate the size and power of test based on $T_{n}$ for this data, and in order to carry out this study, we follow the same methodology as we did for the Canadian weather data in Section \ref{RDA}. Using that methodology, at 5\% level of significance, the estimated size is $0.058$, and the estimated power is $0.587$. It further indicates the proposed test based on $T_{n}$ can achieve the nominal level and perform good in terms of power. 

\begin{figure}
	\begin{center}
	\includegraphics[width= 4.5in, height = 3.25in]{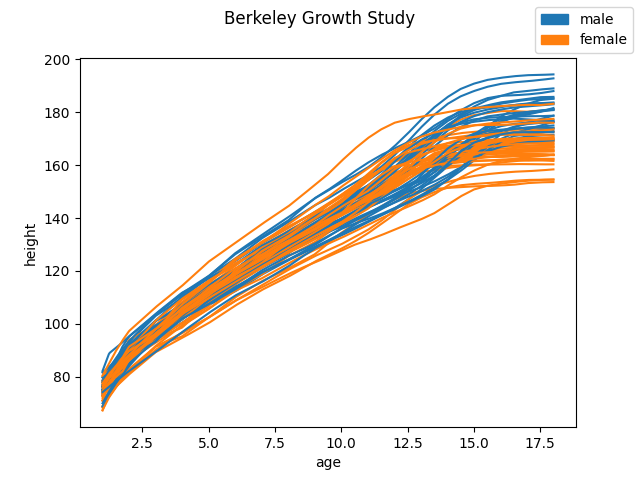}
	\end{center}
	\caption{\it  
		\label{Berkeleygrowth}
		The growth of heights : male and female. Source : \url{https://fda.readthedocs.io/en/latest/auto_examples}} 
\end{figure}

The Berkeley growth data is available at \url{https://rdrr.io/cran/fda/man/growth.html}. It contains the heights of 39 boys and 54 girls measured at 31 time points between the ages 1 and 18 years, and the curves are recorded at 101 equispaced ages in the interval $[1, 18]$. The heights are illustrated in the diagrams in Figure \ref{Berkeleygrowth}, and each observation is viewed as an element in the separable Hilbert space $L_{2}([0, 1])$ after
normalization by an appropriate scaling factor. For this data set, we obtain $p$-value as $0.074$, which indicates that the data
does not favour the null hypothesis at $8\%$ level of significance, i.e., in other words, the heights of boys and girls have strong enough statistical dependence. Note that this result is not unexpected as it seems from the diagrams in Figure \ref{Berkeleygrowth} that the
height curves associated with boys and girls
are not independent. Now, as this data is not favouring the null hypothesis, here also, we
estimate the size and power of the test based on $T_{n}$. At $5\%$ level of significance, the estimated size is $0.057$, and the estimated power is $0.788$, which further indicates that the proposed test is capable of detecting dependence structure of the random elements in real data also. 

\subsection{Simulation Studies}
Here we study the performance of the proposed test for a few more cases when the sample size is finite. The finite sample power of the test is estimated for the following examples. 

\noindent {\bf Example 7 :} Let $(X, Y)$ be $L_{2}([0, 1])\times L_{2}([0, 1])$ valued random element. Suppose that $X(t)$ ($t\in [0, 1]$) follows a Gaussian process $B(t)$ with  $E(B(t))= 0$ and $E(B(t)B(s)) = \min(t, s)$ for all $t\in [0, 1]$ and $s\in [0, 1]$, and $Y(t)\stackrel{d}= \{X(t)\}^{2}$.

\noindent {\bf Example 8 :} Let $(X, Y)$ be $L_{2}([0, 1])\times L_{2}([0, 1])$ valued random element. Suppose that $X(t)$ ($t\in [0, 1]$) follows a Gaussian process $B(t)$ with  $E(B(t))= 0$ and $E(B(t)B(s)) = \min(t, s)$ for all $t\in [0, 1]$ and $s\in [0, 1]$, and $Y(t)\stackrel{d}= e^{X(t)}$.  
 
\noindent {\bf Example 9 :} Let $(X, Y)$ be $L_{2}([0, 1])\times L_{2}([0, 1])$ valued random element. Suppose that $X(t)$ ($t\in [0, 1]$) follows a $t$-process with 3 degrees of freedom, where $t$-process with $k (k\geq 1)$ degrees of freedom is defined as $\frac{B(t)}{\sqrt{\frac{N}{k}}}$, where  $B(t)$ is a Gaussian process  with  $E(B(t))= 0$ and $E(B(t)B(s)) = \min(t, s)$ for all $t\in [0, 1]$ and $s\in [0, 1]$ and independent of $N$, which follows a Chi-squared distribution with $k$ degrees of freedom. Here, $Y(t)\stackrel{d}= \{X(t)\}^{2}$.

\noindent {\bf Example 10 :} Let $(X, Y)$ be an $L_{2}([0, 1])\times L_{2}([0, 1])$ valued random element. Suppose that $X(t)$ ($t\in [0, 1]$) follows a $t$-process with 3 degrees of freedom. Here, $Y(t)\stackrel{d}= e^{X(t)}$.  

For all these aforementioned examples, in order to generate data from the Brownian motion, the data are generated from associated multivariate normal distributions. We here study the performance of the proposed test for the sample sizes $n = 20$, 50, 100 and 500, and the estimated powers of the proposed test are reported in Table \ref{tab3}. 

\begin{table}[h!]
	
	\begin{center}		
		
		\begin{tabular}{ccccc}\hline
			model & $n = 20$ & $n = 50$ & $n = 100$ & $n = 500$\\ \hline 
			Example 7 ($\alpha = 5\%$) & {$0.441$} & {$0.553$} & {$0.685$} & {$0.832$}\\ \hline
			Example 7 ($\alpha = 10\%$) & {$0.483$} & {$0.605$} & {$0.733$} & {$0.886$}\\ \hline
		    Example 8 ($\alpha = 5\%$) & {$0.459$} & {$0.601$} & {$0.824$} & {$0.954$}\\ \hline
		    Example 8 ($\alpha = 10\%$) & {$0.499$} & {$0.643$} & {$0.877$} & {$0.965$}\\ \hline
		    Example 9 ($\alpha = 5\%$) & {$0.484$} & {$0.655$} & {$0.872$} & {$0.963$}\\ \hline
		    Example 9 ($\alpha = 10\%$) & {$0.511$} & {$0.688$} & {$0.901$} & {$0.983$}\\ \hline
		    Example 10 ($\alpha = 5\%$) & {$0.502$} & {$0.713$} & {$0.907$} & {$0.986$}\\ \hline
		    Example 10 ($\alpha = 10\%$) & {$0.532$} & {$0.741$} & {$0.933$} & {$0.992$}\\ \hline
		\end{tabular}
	\end{center}
	\caption{\it The estimated power of the proposed test for different sample sizes $n$. The level of significance, i.e., $\alpha$ are $5\%$ and $10\%$.}
	
	\label{tab3}
	
\end{table}
\noindent The results reported in Table \ref{tab3} indicate that the test is more powerful when $X$ follows $t$-process with 3 degrees of freedom than when $X$ follows Brownian motion. It may be due to the fact that for a fixed $t$, $X(t)$ distributed with a $t$-process with 3 degrees of freedom has a heavier tail than $X(t)$ distributed with a certain Gaussian process. In terms of the relationship between $Y$ and $X$, it is observed that the proposed test is more powerful when $Y \stackrel{d} = e^{X}$ than when $Y\stackrel{d}= X^2$. It is also expected as the exponential function has a faster growth rate than polynomial function.

\noindent{\bf Acknowldgement:}
The first author is thankful to the research grant DST/INSPIRE/04/2017/002835, Government of India and the second author is thankful to the research grants  MTR/2019/000039 and CRG/2022/001489, Government of India. 

\bibliography{lit}

\begin{thebibliography}{}

\bibitem[Bergsma and Dassios, 2014]{Bergsma2014}
Bergsma, W. and Dassios, A. (2014).
\newblock A consistent test of independence based on a sign covariance related
  to kendall’s tau.
\newblock {\em Bernoulli}, 20(2):1006--1028.

\bibitem[Berrett et~al., 2021]{Berrett2021}
Berrett, T.~B., Kontoyiannis, I., and Samworth, R.~J. (2021).
\newblock Optimal rates for independence testing via u-statistic permutation
  tests.
\newblock {\em The Annals of Statistics}, 49(5):2457--2490.

\bibitem[Blum et~al., 1961]{Blum1961}
Blum, J.~R., Kiefer, J., and Rosenblatt, M. (1961).
\newblock Distribution free tests of independence based on the sample
  distribution function.
\newblock {\em The Annals of Mathematical Statistics}, 32:485--498.

\bibitem[Chatterjee, 2021]{Chatterjee2021}
Chatterjee, S. (2021).
\newblock A new coefficient of correlation.
\newblock {\em Journal of the American Statistical Association},
  116(536):2009--2022.

\bibitem[Cuesta-Albertos et~al., 2006]{Juan2006}
Cuesta-Albertos, J.~A., Friman, R., and Ransford, T. (2006).
\newblock Random projections and goodness-of-fit tests in infinite-dimensional
  spaces.
\newblock {\em Bulletin of the Brazilian Mathematical Society}, 37(4):477--501.

\bibitem[Dette et~al., 2013]{Dette2013}
Dette, H., Siburg, K.~F., and Stoimenov, P.~A. (2013).
\newblock A copula-based non-parametric measure of regression dependence.
\newblock {\em Scandinavian Journal of Statistics}, 40(1):21--41.

\bibitem[Dhar et~al., 2018]{Dhar2018}
Dhar, S.~S., Bergsma, W., and Dassios, A. (2018).
\newblock Testing independence of covariates and errors in non-parametric
  regression.
\newblock {\em Scandinavian Journal of Statistics}, 45:421--443.

\bibitem[Dhar et~al., 2016]{DharEJS2016}
Dhar, S.~S., Dassios, A., and Bergsma, W. (2016).
\newblock A study of the power and robustness of a new test for independence
  against contiguous alternatives.
\newblock {\em Electronic Journal of Statistics}, 10:330--351.

\bibitem[Drton et~al., 2020]{Drton2020}
Drton, M., Han, F., and She, H. (2020).
\newblock High-dimensional consistent independence testing with maxima of rank
  correlations.
\newblock {\em The Annals of Statistics}, 48(6):3206--3227.

\bibitem[Duong and Hazelton, 2005]{Duong2005}
Duong, T. and Hazelton, M. (2005).
\newblock Cross validation bandwidth matrices for multivariate kernel density
  estimation.
\newblock {\em Scandinavian Journal of Statistics}, 32:485--506.

\bibitem[Einmahl and Van~Keilegom, 2008]{Einmahl2008}
Einmahl, J. H.~J. and Van~Keilegom, I. (2008).
\newblock Tests for independence in nonparametric regression.
\newblock {\em Statistica Sinica}, 18:601--615.

\bibitem[Ferraty and Vieu, 2006]{Fer2006}
Ferraty, F. and Vieu, P. (2006).
\newblock {\em Nonparametric functional data analysis: Theory and practice}.
\newblock Springer.

\bibitem[Genest et~al., 1961]{Genest2007}
Genest, C., Quessy, J.~F., and Remillard, B. (1961).
\newblock Asymptotic local efficiency of cramer-von mises tests for
  multivariate independence.
\newblock {\em The Annals of Statistics}, 35:166--191.

\bibitem[Gretton et~al., 2007]{Gretton2007}
Gretton, A., Fukumizu, K., Teo, C.~H., Song, L., Scholkopf, B., and Smola,
  A.~J. (2007).
\newblock A kernel statistical test of independence.
\newblock {\em Advances in Neural Information Processing Systems}, pages 1--8.

\bibitem[Hall et~al., 1992]{Hall1992}
Hall, P., Marron, J., and Park, B. (1992).
\newblock Smoothed cross-validation.
\newblock {\em Probability Theory and Related Fields}, 92:1--20.

\bibitem[Han et~al., 2017]{Han2017}
Han, F., Chen, S., and Liu, H. (2017).
\newblock Distribution-free tests of independence in high dimensions.
\newblock {\em Biometrika}, 104(4):813--828.

\bibitem[Lyons, 2013]{Russell2013}
Lyons, R. (2013).
\newblock Distance covariance in metric spaces.
\newblock {\em The Annals of Probability}, 41(5):3284--3305.

\bibitem[Lyons, 2018]{Russell2018}
Lyons, R. (2018).
\newblock Errata to ``distance covariance in metric spaces".
\newblock {\em The Annals of Probability}, 46(4):2400--2405.

\bibitem[Lyons, 2021]{Russell2021}
Lyons, R. (2021).
\newblock Errata to ``distance covariance in metric spaces".
\newblock {\em The Annals of Probability}, 49(5):2668--2670.

\bibitem[Ramsay and Silverman, 2002]{Ramsey2002}
Ramsay, J. and Silverman, B.~W. (2002).
\newblock {\em Applied Functional Data Analysis: Methods and Case Studies}.
\newblock Springer.

\bibitem[Ramsay and Silverman, 2005]{Ramsey2005}
Ramsay, J. and Silverman, B.~W. (2005).
\newblock {\em Functional data analysis}.
\newblock Springer.

\bibitem[Rudin, 1976]{Rudin1976}
Rudin, W. (1976).
\newblock {\em Principles of mathematical analysis}.
\newblock International Series in Pure and Applied Mathematics. McGraw-Hill
  Book Co., New York-Auckland-D\"{u}sseldorf, third edition.

\bibitem[Rudin, 1991]{Rudin1991}
Rudin, W. (1991).
\newblock {\em Functional Analysis}.
\newblock McGraw-Hill.

\bibitem[Serfling, 1980]{Serfling1980}
Serfling, R. (1980).
\newblock {\em Approximation Theorems of Mathematical Statistics}.
\newblock John Wiley \& Sons.

\bibitem[She et~al., 2022a]{Shi2022JASA}
She, H., Drton, M., and Han, F. (2022a).
\newblock Distribution-free consistent independence tests via center-outward
  ranks and signs.
\newblock {\em Journal of the American Statistical Association},
  117(537):395--410.

\bibitem[She et~al., 2022b]{Shi2022}
She, H., Hallin, M., Drton, M., and Han, F. (2022b).
\newblock On universally consistent and fully distribution-free rank tests of
  vector independence.
\newblock {\em The Annals of Statistics}, 50(4):1933--1959.

\bibitem[Silverman, 1998]{Silverman1998}
Silverman, B.~W. (1998).
\newblock {\em Density Estimation for Statistics and Data Analysis}.
\newblock Taylor \& Francis, CRC Press.

\bibitem[Szekely et~al., 2007]{Szekely2007}
Szekely, G.~J., Rizzo, M.~L., and Bakirov, N.~K. (2007).
\newblock Measuring and testing dependence by correlation of distances.
\newblock {\em The Annals of Statistics}, 35(6):2769--2794.

\bibitem[van~der Vaart, 1998]{van1998}
van~der Vaart, A.~W. (1998).
\newblock {\em Asymptotic Statistics}.
\newblock Cambridge University Press.

\end{thebibliography}

\end{document}